\documentclass[reqno,a4paper,11pt]{amsart}

\usepackage{amsmath,amstext,amsfonts,amsbsy,eucal,amssymb}
\usepackage[latin1]{inputenc}
\usepackage[dvips]{epsfig}

\numberwithin{equation}{section}

\newtheorem{theorem}{Theorem}[section]
\newtheorem{lemma}[theorem]{Lemma}

\theoremstyle{definition}
\newtheorem{definition}[theorem]{Definition}
\newtheorem{example}[theorem]{Example}
\newtheorem{remark}[theorem]{Remark}

\begin{document}

\parskip 4pt
\baselineskip 16pt


\title[On the Stability of Some Spline Collocation Implicit Difference Scheme]
{On the Stability of Some Spline Collocation Implicit Difference Scheme}

\author[Svetlana Gaidomak]{Svetlana Gaidomak}
\address{Svetlana Gaidomak, 
Institute for System Dynamics and Control Theory, Siberian Branch of Russian Academy of Sciences,
P.O. Box 292, 664033 Irkutsk, Russia}
\email{gaidamak@icc.ru}
%
%

\date{\today}



\begin{abstract}
Boundary problem for linear partial differential algebraic equations system with multiple characteristic curves is considered. It is supposed that matrix-functions pencil of the system under consideration is smoothly equivalent to special canonical form. For this problem, with the help of the spline collocation method, a difference scheme of arbitrary degree of approximation with respect to each variable is constructed. Sufficient conditions for its absolute stability is found.
\end{abstract}

\maketitle

\section{Introducton}
\label{sect1}

When modeling some processes of hydrodynamics, gas dynamics, atmospheric physics, plasma physics etc., it appears systems with identically degenerate matrix-functions in its domain of definition at all higher partial derivatives \cite{dem}-\cite{han}. Such equations are known in the literature as partial differential algebraic equations, equations not resolved with respect to higher derivatives, degenerate systems of partial differential equations or Sobolev's equations. There exists a number of lines of research of these equations. One of these lines is based on investigation of canonical structures of matrix pencils \cite{cam},~\cite{luch1},~\cite{gai1}-\cite{gai4}. At the moment in the literature, it is well studied the issues of the existence and numerical solution of partial differential algebraic equations with constant matrix coefficients \cite{cam},~\cite{luch},~\cite{han},~\cite{chistbtv}. It is explained by the fact that canonical structures of constant matrix pencils are well studied. Insufficiently studied global properties of matrix-functions pencils impede the investigation of the systems with variable matrix coefficients. In the works \cite{gai1}-\cite{gai4} we considered the systems with simple characteristic curves whose matrix-functions pencils satisfy the criterion known as ``rank-degree'' or double ``rank-degree'' one. Matrices in these pencils depend only on two variables \cite{chist}. For such systems we constructed two-layer and three-layer implicit difference schemes  of a first and second orders of approximation \cite{gai1}-\cite{gai4}. This accuracy is not always enough, for example, in the case when the Lipschitz constant and a domain of definition are quite large. In this case, it requires difference schemes of higher orders of approximation. In the work \cite{gai}, we studied matrix-functions pencils depending on many variables. In a result we obtained sufficient conditions of a smoothly equivalence of these pencils to canonical structure similar to Kronecker form of a regular pencil of constant matrices. This structure is in fact a generalization of a canonical form of the pencil satisfying to the ``rank-degree'' criterion. This structure allows to investigate linear degenerate hyperbolic systems with arbitrary number of independent variables with multiple characteristic curves.

In this paper we consider linear partial differential algebraic equations system of hyperbolic type with the pencil smoothly equivalent to its canonical form similar to the Kronecker form \cite{gai}. With the help of spline collocation method, the foundations of which are set out in \cite{zav}, we construct a high-performance implicit difference scheme of higher order of approximation and then we prove its stability.

The rest of the paper is organized as follows. In the next section, we give a statement of the problem. In section \ref{sect3}, we perform an approximation of unknown function on the uniform grid by a spline of degree $m_{1}$ and $m_{2}$ with respect to each variable, respectively. Then we write down our difference scheme. This approach for constructing of difference scheme was used in \cite{gai4}. In the work \cite{gai4}, for an approximation of unknown function we used the bi-cubic spline with defect. Accuracy of the calculations in this case was low. It is worth remarking that in the present work the corresponding results significantly improved. In section \ref{sect4}, we give some  notations  and auxiliary propositions necessary  to justify correctness of difference scheme are given. The section \ref{sect5} is designed to cast the difference scheme to some canonical form. In this form one can easily analyze the spectrum of matrix coefficients. It is shown that for any values of steps,  the spectrum is entirely contained in  the circle of unit radius. In section \ref{sect6} we prove an absolute  stability property of our difference scheme. Finally, in section \ref{sect7}, for two test examples we show  results of  numerical experiments.

\section{Statement of the problem}
\label{sect2}

Let us consider a boundary problem for the linear partial differential  equations system
\begin{equation}
A(x,t)\partial_{t} u+B(x,t)\partial_{x} u+C(x,t)u=f,
\label{e1}
\end{equation}
\begin{equation}
u(x_{0},t)=\psi(t),\  \ u(x,t_{0})=\phi(x),
\label{e2}
\end{equation}
where  $A(x,t)$, $B(x,t)$ and $C(x,t)$  are  $n\times n$ matrices whose elements are supposed to depend on $x\in \mathbb{R}^{1}$ and  $t\in \mathbb{R}^{1}$. In the following we suppose that $(x,t)\in U=[x_{0};X]\times[t_{0};T]$  and $U$ is therefore a domain of definition of all functions under consideration.
It is supposed that the elements of matrices $A(x,t),\ B(x,t),\ C(x,t)$ and free term $f(x,t)$ belong to $C^{2}(U)$. The vectors $\psi(t)$ and $\phi(x)$ are supposed to be some given vector-functions of its arguments whose elements belong to $C^{2}([t_0, t])$ and $C^{2}([x_0, X])$, respectively.

Let us suppose now that
\begin{equation}
\det A(x,t)=0\ \ \mbox{and} \ \ \det B(x,t)=0\ \ \ \forall \ (x,t)\in U.
\label{ey1}
\end{equation}
The system (\ref{e1}) with the condition (\ref{ey1}) is therefore partial differential algebraic equations one. Its investigation is closely related to analysis of global properties of the pencil $A(x,t)+\lambda B(x,t)$. In the work \cite{gai}, we obtained the conditions of smooth equivalence of matrix-functions pencil to its canonical form similar to canonical structure of regular pencil of constant matrix. in this connection, let us remember the definition of smooth equivalency of  matrix-functions pencils.

\begin{definition} \cite{gai}
Two $n\times n$ matrices pencils $A(x,t)+\lambda B(x,t)$ and $\tilde A(x,t)+\lambda \tilde B(x,t)$ with elements belonging to $C^{s}(U)$, where $\lambda$ is a some parameter are called $s$-smoothly equivalent if there exist square matrices $P(x,t)$ and $Q(x,t)$  which do not depend on $\lambda$ and satisfy following conditions:
\begin{enumerate}
\item the elements of matrix $P(x,t)$ and $Q(x,t)$ belong to $C^{s}(U)$;
\item $\forall (x,t)\in U$ there exist $P^{-1}(x,t)$ and $Q^{-1}(x,t)$;
\item the relation $P(x,t)(A(x,t)+\lambda B(x,t))Q(x,t)=\tilde A(x,t)+\lambda \tilde B(x,t)$ $\forall (x,t)\in U$ holds.
\end{enumerate}
\label{o1}
\end{definition}
It is useful following theorem \cite{gai}.

\begin{theorem} \cite{gai}
Let the following conditions be fulfilled:
\begin{enumerate}
\item all  roots of the characteristic polynomial  $\det (A(x,t)+\lambda B(x,t))$ are real and have a constant multiplicity in a domain of definition $U$;
\item  ranks of matrices $A(x,t)$ and $B(x,t)$ are constant at each point of a domain $U$ and less than $n$.
\end{enumerate}
Then the pencil  $A(x,t)+\lambda\ B(x,t)$ is  smoothly equivalent to the canonical one
\begin{equation}
{\rm diag}\{E_{d},{M}(x,t), E_{p}\}+\lambda\ {\rm diag}\{{J}(x,t), E_{l},{N}(x,t)\},
\label{ee1}
\end{equation}
where  $E_{d}$ is an identity matrix of an order $d$; ${M}(x,t)$ and ${N}(x,t)$ are  an upper (right) triangular  blocks with zero diagonal of orders $l$ and $p$,  respectively; $\mathcal{O}_{l}$ is a zero square matrix of order $l$; $J(x,t)={\rm diag}\{ J_{1}(x,t),  J_{2}(x,t),\ldots,  J_{k}(x,t)\}$, where $J_{i}(x,t)$, for $i=1,\ldots,k$ are nonsingular matrices of orders $d_{i}$, respectively; $d=\sum\limits_{\nu=1}^{k}d_{\nu}$; each block $J_{i}$ has a unique eigenvalue $-1/\lambda_{i}(x,t)$ in the domain of definition $U$; $\lambda_{i}(x,t)$ are eigenvalues of characteristic polynomial  $\det(A(x,t)+\lambda B(x,t))$, different from zero in the domain of definition $U$; $p=n-d-l$.
\label{tp}
\end{theorem}
\begin{remark}
Let conditions of the theorem \ref{tp} are valid. Then if one requires
\[
{\rm rank}(B(x,t))=\deg \left (\det(A(x,t)+\lambda B(x,t))\right )
\]
or
\[
{\rm rank}(A(x,t))=\deg \left (\tilde\lambda\det(A(x,t)+B(x,t))\right ),
\]
then  the pencil $A(x,t)+\lambda B(x,t)$  is $s$-smoothly equivalent to the pencil (\ref{ee1}), in which  $N(x,t)\equiv\mathcal{O}_{l}$ or   $M(x,t)\equiv\mathcal{O}_{p}$, respectively. In this case one says that the pencil $A(x,t)+\lambda B(x,t)$ satisfy the criterion  ``rank-degree''. A structure of such a pencil was investigated  in \cite{chist}.
\label{z1}
\end{remark}
Let us suppose that matrix-functions pencil  $A(x,t)+\lambda B(x,t)$ of the system (\ref{e1}) satisfy the conditions of theorem \ref{tp} and remark \ref{z1}. In the next section, we proceed to construct our difference scheme.

\section{Difference scheme}
\label{sect3}

To construct the difference scheme , we perform the partition  of the domain $U$ by the lines
$x=x_{i}$, where $x_{i}=x_{0}+ih$ and $t=t_{j}$, where $t_{j}=t_{0}+j\tau$ for $i=0,\ldots, n_{1}$ and $j=0,\ldots, n_{2}$.
In a result, we obtain a uniform grid $U_{\Delta}$ with the steps $h$ and $\tau$ with respect to the space and time variable, respectively. Clearly,
$h=(X-x_{0})/{n}_{1}$ and $\tau=(T-t_{0})/{n}_{2}$, where $0<\tau\leq \tau_{0}$ and $0<h\leq h_{0}$. The points with coordinates $(x_{i},t_{j})$ are referred  to as the nodes of the grid  $U_{\Delta}$ while the lines $x=x_{i}$ and $t=t_{j}$ are called its layers. In each domain $U_{i,j}=[x_{i}, x_{i}+m_{1}h]\times [t_{j}, t_{j}+m_{2}\tau]\subseteq U_{\Delta}$, where $m_{1}\leq {n}_{1}$ and $m_{2}\leq {n}_{2}$, we seek the approximation of the solution $u(x,t)$ of the problem (\ref{e1}),~(\ref{e2}) in the form of a polynomial $L_{i,j}^{m_{1},m_{2}}(x,t)$ with indeterminate coefficients
of degree $m_{1}$  and $m_{2}$ with respect to variables $x$ and $t$, respectively. We require that the values of this polynomial $L_{i,j}^{m_{1},m_{2}}(x,t)$ in the nodes $(x_{i}+l_{1}h, t_{j}+l_{2}\tau)$ for $l_{1}=0,\ldots, m_{1}$ and $l_{2}=0,\ldots, m_{2}$ of a grid domain $U_{i,j}$ coincide with the values of the desired solution $u(x,t)$ at these nodes. In order our approximation will be continuous, we require that the relations
\[
L_{i,j}^{m_{1},m_{2}}(x_{i}+l_{1}h,t_{j})=L_{i+1,j}^{m_{1},m_{2}}(x_{i}+l_{1}h,t_{j}),\ \ L_{0,j}^{m_{1},m_{2}}(x_{0},t_{j})=\psi(t_{j}),
\]
\[
L_{i,j}^{m_{1},m_{2}}(x_{i},t_{j}+l_{2}\tau)=L_{i,j+1}^{m_{1},m_{2}}(x_{i},t_{j}+l_{2}\tau),\ \ L_{i,0}^{m_{1},m_{2}}(x_{i},t_{0})=\phi(x_{i}).
\]
hold  at horizontal and vertical line $x=x_{i}$ and $t=t_{j}$, respectively. Applying differenceless formulas of numerical differentiation for equidistant nodes (\cite{ber}, p. 161), to derivatives $\partial_{t} u(x,t)$ and $\partial_{x} u(x,t)$ at layers $x=x_{i}$ and $t=t_{j}$, we obtain
\begin{equation}
\partial_{x}u(x_{i}+l_{1}h,t_{j})=\partial_{x}L_{i,j}^{m_{1},m_{2}}(x_{i}+l_{1}h,t_{j})+\epsilon_{1}(h^{m_{1}}),\ \ l_{1}=1,\ldots, m_{1},
\label{e3}
\end{equation}
\begin{equation}
\partial_{t}u(x_{i},t_{j}+l_{2}\tau)=\partial_{t}L_{i,j}^{m_{1},m_{2}}(x_{i},t_{j}+l_{2}\tau)+\epsilon_{2}(\tau^{m_{2}}),\ \ l_{2}=1,\ldots, m_{2},
\label{e4}
\end{equation}
where
\[
\partial_{x}L_{i,j}^{m_{1},m_{2}}(x_{i}+l_{1}h,t_{j})=\frac{1}{h}\sum_{l_{3}=0}^{m_{1}}\bar\gamma_{l_{1},l_{3}}u(x_{i}+l_{1}h,t_{j}),
\]
\[
\partial_{t}L_{i,j}^{m_{1},m_{2}}(x_{i},t_{j}+l_{2}\tau)=\frac{1}{\tau}\sum_{l_{3}=0}^{m_{2}}\gamma_{l_{2},l_{3}}u(x_{i},t_{j}+l_{3}\tau),
\]
\[
\epsilon_{1}(h^{m_{1}})=\frac{h^{m_{1}}{\partial^{m_{1}+1}(\zeta_{1},t_{j})}/{\partial x^{m_{1}+1}}}{(m_{1}+1)!}\left.\frac{d}{ds}\prod_{\nu=0}^{m_{1}}(s-\nu)\right|_{s=l_{1}},\ \ x_{i}<\zeta_{1}<x_{i+1},
\]
\[
\epsilon_{2}(\tau^{m_{1}})=\frac{\tau^{m_{2}}{\partial^{m_{2}+1}(x_{i},\zeta_{2})}/{\partial t^{m_{2}+1}}}{(m_{2}+1)!}\left.\frac{d}{ds}\prod_{\nu=0}^{m_{2}}(s-\nu)\right|_{s=l_{2}},\ \ t_{j}<\zeta_{2}<t_{j+1}.
\]
The coefficients $\bar\gamma_{l_{1},l_{3}}$  and $\gamma_{l_{2},l_{3}}$  are calculated with the help of following formulas
 (\cite{ber}, p. 161):
\[
\bar\gamma_{l_{1},l_{3}}=\left.{ H}(m,s,l_{3})\right |_{m=m_{1},\ s=l_{1}},\ \
\gamma_{l_{2},l_{3}}=\left.{ H}(m,s,l_{3})\right |_{m=m_{2},\ s=l_{2}},
\]
\begin{equation}
{ H}(m,s,l_{3})=(-1)^{m+l_{3}}\frac{C^{l_{3}}_{m}}{m!}\frac{d}{d s} \left ({ \prod\limits_{\nu=0}^{m}(s-\nu)}/{(s-l_{3})}\right ),
\label{e5}
\end{equation}
where $C^{l_{3}}_{m}$ stand for the binomial coefficients. Writing the system (\ref{e1}) at the nodes $(x_{i}+l_{1}h,t_{j}+l_{2}\tau)$ for $l_{1}=1,\ldots, m_{1}$ and $l_{2}=1,\ldots, m_{2}$ of the domain $U_{i,j}$ and substituting into it the values of desired function $u(x_{i}+l_{1}h,t_{j}+l_{2}\tau)$ and approximation of its derivatives (\ref{e3}) and (\ref{e4}) in these points, we obtain following difference scheme:
\[
A_{i+l_{1},j+l_{2}}\frac{1}{\tau}\sum_{l_{3}=0}^{m_{2}}\gamma_{l_{2},l_{3}}u_{i+l_{1},j+l_{3}}+B_{i+l_{1},j+l_{2}}\frac{1}{h}\sum_{l_{3}=0}^{m_{1}}\bar\gamma_{l_{1},l_{3}}u_{i+l_{3},j+l_{2}}
\]
\[
+C_{i+l_{1},j+l_{2}}u_{i+l_{1},j+l_{2}}=f_{i+l_{1},j+l_{2}},\\
\]
\begin{equation}
u_{0,j}=\psi_{j},\;\; u_{i,0}=\phi_{i},\;\;  i=0,\ldots,n_{1}-1,\;\;  i=0,\ldots,n_{2}-1,
\label{e6}
\end{equation}
where
\[
A_{i+l_{1},j+l_{2}}\equiv A(x_{i}+l_{1}h,t_{j}+l_{2}\tau),\;\;
B_{i+l_{1},j+l_{2}}\equiv B(x_{i}+l_{1}h,t_{j}+l_{2}\tau),\;\;
\]
\[
C_{i+l_{1},j+l_{2}}\equiv C(x_{i}+l_{1}h,t_{j}+l_{2}\tau),\;\;
f_{i+l_{1},j+l_{2}}\equiv f(x_{i}+l_{1}h,t_{j}+l_{2}\tau),\;\;
\]
\[
u_{i+l_{1},j+l_{2}}\equiv u(x_{i}+l_{1}h,t_{j}+l_{2}\tau),\;\;
\phi_{i}\equiv\phi(x_{i}), \psi_{j}=\psi(t_{j}).
\]
The difference scheme (\ref{e6}) at each node of the grid $U_{\Delta}$ is given by a system of linear algebraic equations of order $\tilde n=m_{1}m_{2}n$ with unknown vector
\[
\bar u_{i+1,j+1}=(u_{i+1,j+1},\dots, u_{i+1,j+m_{2}},u_{i+2,j+1}, \dots, u_{i+2,j+m_{2}},\dots ,u_{i+m_{1},j+1}, \dots, u_{i+m_{1},j+m_{2}})^{\top}.
\]
In what follows to avoid a confusion in notations, we denote the solution of the system (\ref{e6}) in the nodes $(x_{i+1},t_{j+1})$ of the grid $ U_{\Delta}$ by $v_{i+1,j+1}$, while the values of unknown function $u(x,t)$ in the same nodes we denote as $u(x_{i+1},t_{j+1})=u_{i+1,j+1}$. Then the system (\ref{e6}) has the following form:
\[
A_{i+l_{1},j+l_{2}}\frac{1}{\tau}\sum_{l_{3}=1}^{m_{2}}\gamma_{l_{2},l_{3}}v_{i+l_{1},j+l_{3}}+B_{i+l_{1},j+l_{2}}\frac{1}{h}\sum_{l_{3}=1}^{m_{1}}\bar\gamma_{l_{1},l_{3}}v_{i+l_{3},j+l_{2}}+C_{i+l_{1},j+l_{2}}v_{i+l_{1},j+l_{2}}=
\]
\begin{equation}
=f_{i+l_{1},j+l_{2}}-\frac{1}{\tau}A_{i+l_{1},j+l_{2}}\gamma_{l_{2},0}v_{i+l_{1},j}-\frac{1}{h}B_{i+l_{1},j+l_{2}}\bar\gamma_{l_{1},0}v_{i,j+l_{2}},\\
\label{e7}
\end{equation}
\[
v_{0,j}=\psi_{j},\;\;
v_{i,0}=\phi_{i},\;\;
i=0,\ldots, n_{1}-1,\;\;
i=0,\ldots, n_{2}-1,\;\;
l_{1}=1,\ldots, m_{1},\;\;
l_{2}=1,\ldots, m_{2},\;\;
\]
where
\[
\bar v_{i+1,j+1}=(v_{i+1,j+1},\ldots, v_{i+1,j+m_{2}},v_{i+2,j+1}, \ldots, v_{i+2,j+m_{2}},\ldots ,v_{i+m_{1},j+1}, \ldots, v_{i+m_{1},j+m_{2}})^{\top}.
\]
Let us remark that the difference scheme (\ref{e7}) in fact represents the whole set of implicit difference schemes. Setting different values of the orders $m_{1}$ and $m_{2}$ of approximating polynomials $L_{i,j}^{m_{1},m_{2}}(x,t)$, we get according to (\ref{e3}) and (\ref{e4}) difference schemes with corresponding approximation orders $O(h^{m_{1}})+O(\tau^{m_{2}})$. In addition, the scheme (\ref{e7})  can be both two-layer and multi-layer depending on given information.

In our case we have in mind two-layer difference scheme, that is, we suppose that the values of unknown function are given only at one left layer $x=x_{i}$ and on one lower layer $t=t_{j}$ of $U_{\Delta}$ (at the nodes of these layers  the values  of grid function were defined in the previous step), but in actual calculations  are used $m_{1}+1$ vertical and $m_{2}+1$ horizontal layers, as is shown on the template (Fig. 1). The movement on the grid is carried along its layers  and within each layer by steps.
\begin{figure}[htb]
\begin{center}
\includegraphics[width=0.45 \textwidth, height=0.3\textwidth ]{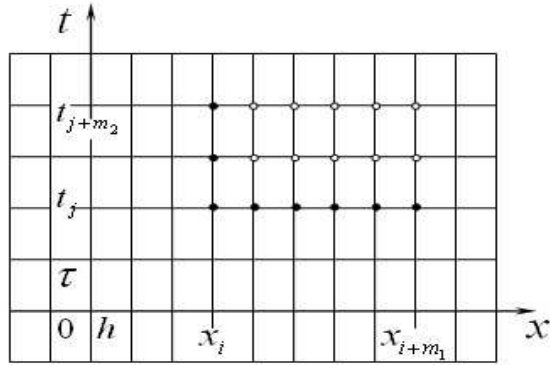}
{\caption{\footnotesize {The grid $U_{\Delta}$ and $(m_{1}+1)(m_{2}+1)$-point template.}}}
\end{center}
\end{figure}
Therefore we have constructed the difference scheme (\ref{e7}) or more exactly  a whole set of difference schemes with different orders of approximation. To proceed, we must to find out the conditions of solvability for the difference scheme (\ref{e7}) for any data and for sufficiently small grid steps and finally to prove its stability following the definition from (\cite{sam}, p. 12). It is worth remarking that the difference scheme (\ref{e7}) was written  in \cite{gai5} but without proof of its stability.

To estimate a solution of the problem (\ref{e7}), in the space $C(U_{\Delta})$ consisting of $n$-dimensional vector-functions $v(x_{i},t_{j})=(v^{1}(x_{i},t_{j}),\ v^{2}(x_{i},t_{j}),\dots,v^{n}(x_{i},t_{j}))^{\top}$  let us use the following norm:
\[
||v(x_{i},t_{j})||_{C(U_{\Delta})}=\max_{(x_{i},t_{j})\in U_{\Delta}}||v(x_{i},t_{j})||,\ \ \mbox{where}\ \ ||v(x_{i},t_{j})||=\max_{\nu=1,\ldots, n} |v^{\nu}(x_{i},t_{j})|.
\]
It is obvious that it is consisted with the norm of  $n$-dimensional vector-function $v(x,t)\in C(U)$  \cite{gant}
\[
||v(x,t)||_{C(U)}=\max\{||v(x,t)||\ \forall (x,t)\in U\}.
\]
In the following, we need in some auxiliary propositions which we formulate and prove in the next section.

\section{Preliminaries}
\label{sect4}

Let us spend some lines to fix auxiliary notations used throughout the paper. Denote $\tau/h=r$ and define diagonal matrices
\[
\bar\gamma^{0}_{m_{1}}={\rm diag}\{\bar\gamma_{1,0},\bar\gamma_{2,0},\dots,\bar\gamma_{m_{1},0}\},\ \ \gamma^{0}_{m_{2}}={\rm diag}\{\gamma_{1,0},\gamma_{2,0},\dots,\gamma_{m_{2},0}\},
\]
\begin{equation}
\bar\gamma_{m_{1}}=(\bar\gamma_{i,j})_{1}^{m_{1}},\ \ \ \gamma_{m_{2}}=(\gamma_{i,j})_{1}^{m_{2}}
\label{em9}
\end{equation}
whose elements are given by (\ref{e5}). Remark, that $\bar\gamma_{m_{1}}$ and $\gamma_{m_{2}}$ coincides if $m_{1}=m_{2}$. Also we define following block diagonal matrices:
\[
\Omega_{i+1,j+1}={\rm diag}\{\Omega_{i+1,j+1}^{1},\ \Omega_{i+1,j+1}^{2},\ \Omega_{i+1,j+1}^{3}\},
\]
\begin{equation}
\bar{\\mathcal{A}}={\rm diag}\{\mathcal{A}^{1}, \mathcal{A}^{2}, \mathcal{A}^{3}\},\ \ \bar{\mathcal{B}}_{i+1,j+1}={\rm diag}\{\mathcal{B}_{i+1,j+1}^{1}, \mathcal{B}^{2}, \mathcal{B}^{3}\},
\label{e9}
\end{equation}
where
\[
\Omega_{i+1,j+1}^{1}=E_{m_{1}}\otimes\gamma_{m_{2}}\otimes E_{d}+r \bar\gamma_{m_{1}}\otimes E_{m_{2}}\otimes {J}_{i+1,j+1},
\]
\[
\Omega^{2}=r \bar\gamma_{m_{1}}\otimes E_{m_{2}l},\ \
\Omega^{3}=E_{m_{1}}\otimes\gamma_{m_{2}}\otimes E_{p}
\]
\[
\mathcal{A}^{1}=E_{m_{1}}\otimes \gamma^{0}_{m_{2}}\otimes E_{d},\ \ \mathcal{A}^{2}=\mathcal{O}_{m_{1}m_{2}l},\ \ \mathcal{A}^{3}=E_{m_{1}}\otimes \gamma^{0}_{m_{2}}\otimes E_{p},
\]
\[
\mathcal{B}_{i+1,j+1}^{1}=\bar \gamma^{0}_{m_{1}}\otimes E_{m_{2}}\otimes {J}_{i+1,j+1},\ \ \mathcal{B}^{2}=\bar\gamma_{m_{1}}^{0}\otimes E_{m_{2}}\otimes E_{l},
\ \ \mathcal{B}^{3}=\mathcal{O}_{m_{1}m_{2}p},
\]
\begin{lemma}
Let $\xi_{\bar\gamma_{m_{1}}}^{s_{1}}$,  $\xi_{\gamma_{m_{2}}}^{s_{2}}$ and $\xi_{{J}_{i+1,j+1}}^{s_{3}}$ for  $s_{1}=1,\ldots, m_{1},$ $s_{2}=1,\ldots, m_{2}$ and $s_{3}=1,\ldots, k$ are eigenvalues of the  matrices $\bar\gamma_{m_{1}}$, $\gamma_{m_{2}}$ and ${J}_{i+1,j+1}$, respectively. If at each node of the grid  $U_{\Delta}$ the inequalities
\begin{equation}
r\xi_{\bar\gamma_{m_{1}}}^{s_{1}}\xi_{{J}_{i+1,j+1}}^{s_{3}}\neq -\xi_{\gamma_{m_{2}}}^{s_{2}},\;\;\; \forall s_{1} ,s_{2}\;\mbox{and}\;s_{3},
\label{e11}
\end{equation}
be fulfilled, then the matrix $\Omega_{i+1,j+1}$  in (\ref{e9}) is nonsingular	on the grid $U_{\Delta}$.
\label{l2}
\end{lemma}

{\sc Proof.}
To prove the nonsingularity of the matrix $\Omega_{i+1,j+1}$ on $U_{\Delta}$ under condition of this lemma, it is sufficient to prove this property for its blocks. Since the matrices $\bar\gamma_{m_{1}}$ and $\gamma_{m_{2}}$ are nonsingular then the blocks $\Omega^{2}$ and $\Omega^{3}$ also share this property. Let us consider now the block $\Omega^{1}_{i+1,j+1}$. Let $\mathcal{P}_1$ and $\mathcal{P}_2$ are constant matrices and $\mathcal{ P}_{3}(x,t)$ is a matrix-function, such that
\[
\mathcal{P}_{1}\bar\gamma_{m_{1}}\mathcal{P}_{1}^{-1}=\tilde\gamma_{m_{1}},\ \
\mathcal{P}_{2}\gamma_{m_{2}}\mathcal{P}_{2}^{-1}=\tilde\gamma_{m_{2}},\ \
\mathcal{P}_{3}(x,t){J}(x,t)\mathcal{P}_{3}^{-1}(x,t)=\tilde{ J}(x,t),
\]
where $\tilde\gamma_{m_{1}}$, $\tilde\gamma_{m_{2}}$ and $\tilde{ J}(x,t)$  are normal Jordan forms of the matrices $\bar\gamma_{m_{1}}$, $\gamma_{m_{2}}$ and ${ J}(x,t)$, respectively. Remark that an existence of the matrix $\mathcal{P}_{3}(x,t)$ is provided by the theorem in \cite{verb}.
Make up $\mathcal{P}(x,t)=\mathcal{P}_{1}\otimes \mathcal{P}_{2} \otimes \mathcal{P}_{3}(x,t)$ and multiply the matrix $\Omega^{1}_{i+1,j+1}$ on the right  and on the left by matrices $\mathcal{P}_{i+1,j+1}$ and $\mathcal{P}_{i+1,j+1}^{-1}$, respectively. In a result, we obtain a normal Jordan form of the block $\Omega^{1}_{i+1,j+1}$. It is obvious, that provided (\ref{e11}) there are no zeros at the main diagonals of Jordan cells of normal form of the matrix  $\Omega^{1}_{i+1,j+1}$. Therefore the lemma is proved.

The following lemma is particularly important because it gives characteristic properties of splines and allows to transform the difference scheme (\ref{e7}) to its canonical form in which  one easily sees spectra of matrix coefficients.
\begin{lemma}
Let $x_{0}\in{\mathbb R}^{d}$ be an arbitrary vector and  $y_{0}=(x_{0},x_{0},\dots, x_{0})^{\top}\in{\mathbb R}^{md}$. Let $\alpha$ be an arbitrary parameter and $\mathcal{J}$ be an arbitrary constant $d\times d$ matrix with corresponding eigenvalues $\xi_{\mathcal{J}}^{l}$ for $l=1,\ldots, d$  satisfying
$\xi_{\mathcal{J}}^{l}\neq -\xi_{\gamma_{m}}^{s}/\alpha,\;\; \forall\ l,s$, where $\xi_{\gamma_{m}}^{s}$ for $s=1,\ldots, m$ are eigenvalues of $m\times m$ matrix $\gamma_m\equiv(\gamma_{i,j})_{1}^{m}$ with elements defined through (\ref{e5}). Then
\begin{equation}
\left ( E_{md}+\alpha (\gamma_{m}^{-1}\otimes \mathcal{J})\right )^{-1}y_{0}=- {\rm diag} \left\{\exp(-\alpha \mathcal{J}), \exp(-2\alpha \mathcal{J}),\dots,\exp(-m\alpha \mathcal{J})\right\}y_{0}+O(\tau^{m}).
\label{e12}
\end{equation}
\label{l3}
\end{lemma}
{\sc Proof.}
Let us consider Cauchy problem for homogeneous system of ordinary differential equations
\begin{equation}
\dot x(t)=-\frac{\alpha}{\tau}\mathcal{J}x(t),\ \ x(t_{0})=x_{0},\ \ t\in I=[0,m\tau].
\label{e13}
\end{equation}
The solution of (\ref{e13}), as is known, is
\begin{equation}
x(t)=\exp\left (-\frac{\alpha}{\tau}\mathcal{J}t\right )x_{0}.
\label{e14}
\end{equation}
Divide the segment $I$ into $(m-1)$ equal parts by the points $t_{j}=j\tau$ for $j=0,\ldots, m$. Making  use of differenceless formulas of numerical differentiation (\ref{e4}), we write down the approximation for the derivative of the vector-function  $x(t)$ at the points  $t_{j}$ for $j=1,\ldots, m$
\begin{equation}
\tau \dot x(j\tau)=\sum_{l=0}^{m}\gamma_{j,l}x(t_{l})+O(\tau^{m}).
\label{e15}
\end{equation}
Rewrite the relations (\ref{e15}) in the matrix form
\begin{equation}
Y=\left(\gamma_{m}^{0}\otimes E_{d}\right)y_{0}+\left(\gamma_{m}\otimes E_{d}\right)y+O(\tau^{m+1}),
\label{e16}
\end{equation}
where $Y$ and $y$ are vectors of the form $Y=(\dot x(h), \dot x(2h), \dots, \dot x(mh))^{\top}$ and $y=(x(h), x(2h), \dots, x(mh))^{\top}$, respectively.
With the help of (\ref{e13}) we find the values of the derivative of unknown function at the points  $t_{j}$  of the segment $I$.
Substituting these values into (\ref{e16}), we get
\begin{equation}
-\alpha (E_{m}\otimes \mathcal{J})y=\left(\gamma_{m}^{0}\otimes E_{d}\right)y_{0}+\left(\gamma_{m}\otimes E_{d}\right)y+O(\tau^{m+1}).
\label{e17}
\end{equation}
Resolving (\ref{e17}) in favor of $y$ gives
\begin{equation}
y=-\left(\alpha (E_{m}\otimes \mathcal{J})+\gamma_{m}\otimes E_{d}\right)^{-1}\left(\gamma_{m}^{0}\otimes E_{d}\right)y_{0}+O(\tau^{m+1}).
\label{e18}
\end{equation}
In the work \cite{gai6}, we proved the identity
\begin{equation}
\left(\gamma_{m}^{-1}\gamma_{m}^{0}\otimes E_{d}\right)e_{md}=-e_{md}+O(\tau^{m}).
\label{el18}
\end{equation}
Here $e_{md}$ stands for a $md$-dimensional vector, each element of which ia a unit. Remark, that relation (\ref{el18}) is in fact a particular case of the following one:
\begin{equation}
\left(\gamma_{m}^{-1}\gamma_{m}^{0}\otimes E_{d}\right)y_{0}=-y_{0}+O(\tau^{m}).
\label{el19}
\end{equation}
Remark, that it is proved in the same way as  (\ref{el18}), but instead of the interpolated vector-function $g(t)=(1+t^{m},1+t^{m},\dots,1+t^{m})^{\top}$,  we take here $\tilde g(t)=x_{0}-e_{d}+g(t)$. Making use of (\ref{el19}), we transform (\ref{e18}) to get
\begin{equation}
y=-\left(E_{md}+\alpha (\gamma_{m}^{-1}\otimes \mathcal{J})\right)^{-1}y_{0}+O(\tau^{m}).
\label{e19}
\end{equation}
Finally, substituting into (\ref{e19}) the values of unknown function (\ref{e14}) at the nodal points of $I$, we obtain the relation (\ref{e12}).  Therefore the lemma is proved.

\section{Transformation of difference scheme under consideration}
\label{sect5}
It is quite difficult to investigate the difference scheme (\ref{e7}) without preliminary transformation of it to convenient form. Thus, let us first to cast it to special canonical form. As a result we write down this canonical form in the end of this section.

We suppose that the pencil $A(x,t)+\lambda B(x,t)$ of the system (\ref{e1}) satisfies the conditions of the theorem \ref{tp}. Remember that this means that one can find a pair of matrix-functions  $P(x,t)$ and $Q(x,t)$ with properties of definition \ref{o1}, which serve for transformation of our pencil  to the canonical form (\ref{ee1}). Prepare following square matrix  of order $\tilde n$:
\[
\tilde P={\rm diag}\left\{P_{i+1,j+1},\dots , P_{i+1,j+m_{2}},P_{i+2,j+1},\dots , P_{i+2,j+m_{2}},\dots ,P_{i+m_{1},j+1},\dots , P_{i+m_{1},j+m_{2}}\right\}.
\]
Multiply the left- and right-hand sides of the system (\ref{e7}) on the left by the matrix $\tau\tilde P$ and perform a change of variable  $v_{i,j}=Q_{i,j}w_{i,j}$, where $w_{i,j}$ is some unknown  $n$-dimensional  vector-function, calculated in the node $(i,j)$. In a result we obtain the following difference scheme:
\[
P_{i+l_{1},j+l_{2}}\Biggl \{ A_{i+l_{1},j+l_{2}}\sum_{l_{3}=1}^{m_{2}}\gamma_{l_{2},l_{3}}Q_{i+l_{1},j+l_{3}}w_{i+l_{1},j+l_{3}}+r B_{i+l_{1},j+l_{2}}\sum_{l_{3}=1}^{m_{1}}\bar\gamma_{l_{1},l_{3}}Q_{i+l_{3},j+l_{2}}w_{i+l_{3},j+l_{2}}
\]
\[
+\tau C_{i+l_{1},j+l_{2}}Q_{i+l_{1},j+l_{2}}w_{i+l_{1},j+l_{2}}\Biggr \}
=
P_{i+l_{1},j+l_{2}}\Biggl \{\tau f_{i+l_{1},j+l_{2}}-rA_{i+l_{1},j+l_{2}}\gamma_{l_{2},0}Q_{i+l_{1},j}w_{i+l_{1},j}
\]
\begin{equation}
-
rB_{i+l_{1},j+l_{2}}\bar\gamma_{l_{1},0}Q_{i,j+l_{2}}w_{i,j+l_{2}}\Biggr \},\\
\label{es7}
\end{equation}
\[
w_{0,j}=Q^{-1}_{0,j}\psi_{j},\ \ w_{i,0}=Q^{-1}_{i,0}\phi_{i},
\]
for $i=0,\ldots, n_{1}-1$, $i=0,\ldots, n_{2}-1$, $l_{1}=1,\ldots, m_{1}$ and $l_{2}=1,\ldots,m_{2}$ with unknown vector
\[
\bar w_{i+1,j+1}=(w_{i+1,j+1},\dots, w_{i+1,j+m_{2}},w_{i+2,j+1}, \dots, w_{i+2,j+m_{2}},\dots ,w_{i+m_{1},j+1}, \dots, w_{i+m_{1},j+m_{2}})^{\top}.
\]
With Taylor's formula we can represent the matrix $Q_{i+l_{1},j+l_{2}}$ in the following form:
\[
Q_{i+l_{1},j+l_{2}}=Q_{i+\tilde l_{1},j+\tilde l_{2}}+\sigma^{1}_{i+l_{1},j+l_{2}}\chi_{1}h+\sigma^{2}_{i+l_{1},j+l_{2}}\chi_{2}\tau,
\]
for $\tilde l_{1}=1,\ldots, m_{1}$ and $\tilde l_{2}=1\ldots ,m_{2}$, where
\[
\sigma^{1}_{i+l_{1},j+l_{2}}\equiv\partial_{x} Q(x_{i+l_{1}}+\theta h,t_{j+l_{2}})\;\;\; \mbox{and}\;\;\; \sigma^{1}_{i+l_{1},j+l_{2}}\equiv\partial_{\tau} Q(x_{i+l_{1}},t_{j+l_{2}}+\theta \tau).
\]
Here $\theta$ is some parameter, by assumption, obeying  the condition $0<\theta<1$ and $\tilde l_{k}$ are suitable numbers for which one has  $\chi_{k}=|\tilde l_{k}-l_{k}|$, while $0\leq \chi_{k}\leq m_{k}$. Let us denote $\tilde C_{i,j}=P_{i,j}C_{i,j}Q_{i,j}$ and rewrite  system (\ref{es7}), taking into account canonical form (\ref{ee1}) of the pencil $A(x,t)+\lambda B(x,t)$. In a result, we obtain
\begin{equation}
\tilde \Omega_{i+1,j+1}\bar w_{i+1,j+1}=q_{i+1,j+1},
\label{es8}
\end{equation}
\[
\tilde\Omega_{i+1,j+1}=\hat\Omega_{i+1,j+1}+h\sigma_{1}+\tau \sigma_{2},
\]
\[
q_{i+1,j+1}=\tau\tilde f_{i+1,j+1}-(\hat{\mathcal{A}}+\tau \sigma_{3})\bar w_{i+1,j}-r(\hat{\mathcal{B}}_{i+1,j+1}+\tau \sigma_{4})\bar w_{i,j+1},
\]
where $\hat\Omega_{i+1,j+1},\ \hat{\mathcal{A}}$ and $\hat{\mathcal{B}}_{i+1,j+1}$ are $\tilde n\times\tilde n$  matrices of the form
\[
\hat\Omega_{i+1,j+1}=E_{m_{1}}\otimes\gamma_{m_{2}}\otimes{\rm diag}\{E_{d}, \mathcal{O}_{l}, E_{p}\}+r \bar\gamma_{m_{1}}\otimes E_{m_{2}}\otimes{\rm diag}\{J_{i+1,j+1}, E_{l}, \mathcal{O}_{p}\},
\]
\[
\hat {\mathcal{A}}=E_{m_{1}}\otimes\gamma_{m_{2}}^{0}\otimes{\rm diag}\{E_{d}, \mathcal{O}_{l}, E_{p}\},\ \
\hat{\mathcal{B}}_{i+1,j+1}=\bar\gamma_{m_{1}}\otimes E_{m_{2}}\otimes{\rm diag}\{J_{i+1,j+1}, E_{l}, \mathcal{O}_{p}\};
\]
and
\[
\tilde f_{i+1,j+1}=\tau \tilde P\bar f_{i+1,j+1},
\]
\[
\bar f_{i+1,j+1}=(f_{i+1,j+1},\dots, f_{i+1,j+m_{2}},f_{i+2,j+1}, \dots, f_{i+2,j+m_{2}},\dots ,f_{i+m_{1},j+1}, \dots, f_{i+m_{1},j+m_{2}})^{\top};
\]
\[
\bar w_{i+1,j}=(e_{m_{2}}\otimes w_{i+1,j},e_{m_{2}}\otimes w_{i+2,j}, \dots, e_{m_{2}}\otimes w_{i+m_{1},j})^{\top},
\]
\[
\bar w_{i,j+1}=e_{m_{1}}\otimes (w_{i,j+1},w_{i,j+2}, \dots, w_{i,j+m_{2}})^{\top}.
\]
Here $\sigma_{k}$ are bounded in the domain of definition matrices obtained as a result of expanding $Q_{i+l_{1},j+l_{2}}$ and ${J}_{i+l_{1},j+l_{2}}$  by Taylor's formula.

Let us split each block component $w_{i+1,j+1}$ of the vector $\bar w_{i+1,j+1}$ into three ones
\[
w_{i+1,j+1}=(w_{i+1,j+1}^{1},w_{i+1,j+1}^{2},w_{i+1,j+1}^{3})^{\top}
\]
of the size $d, l$ and $p$, respectively. Prepare then following $\tilde n\times\tilde n$ permutation matrix
\[
T={\rm colon}(T_{1},T_{2},T_{3}),
\]
\begin{equation}
T_{1}=E_{m_{1}}\otimes E_{m_{2}}\otimes (E_{d}\ \mathcal{O}_{d\times l}\ \mathcal{O}_{d\times p}),\ \ \ T_{2}=E_{m_{1}}\otimes E_{m_{2}}\otimes (\mathcal{O}_{l\times d}\ E_{l}\ \mathcal{O}_{l\times p}),
\label{e20}
\end{equation}
\[
T_{3}=E_{m_{1}}\otimes E_{m_{2}}\otimes (\mathcal{O}_{p\times d}\ \mathcal{O}_{p\times l}\ E_{p}),
\]
where the blocks  $T_{1},\ T_{2}$ and $T_{3}$  are supposed to have the sizes $m_{1}m_{2}d\times \tilde n$,  $m_{1}m_{2}l\times \tilde n$ and $m_{1}m_{2}p\times \tilde n$, respectively. It is easy to prove that $T$ is orthogonal matrix. The matrix $T$ perform the permutation of the components of the vector $\bar w_{i+1,j+1}$ so that
\[
z_{i+1,j+1}\equiv T\bar w_{i+1,j+1}=
(\bar w_{i+1,j+1}^{1},\bar w_{i+1,j+1}^{2},\bar w_{i+1,j+1}^{3})^{\top},
\]
where
\[
\bar w_{i+1,j+1}^{k}=(w_{i+1,j+1}^{k},\dots, w_{i+1,j+m_{2}}^{k},w_{i+2,j+1}^{k}, \dots, w_{i+2,j+m_{2}}^{k},\dots ,w_{i+m_{1},j+1}^{k}, \dots, w_{i+m_{1},j+m_{2}}^{k})^{\top}.
\]
One sees that the vectors  $\bar w_{i+1,j+1}^{1}$, $\bar w_{i+1,j+1}^{2}$ and $\bar w_{i+1,j+1}^{3}$ have
the sizes $m_{1}m_{2}d$, $m_{1}m_{2}l$ and $m_{1}m_{2}p$, respectively.
With the help of  $T$ we also reshuffle the elements of  $\hat \Omega_{i+1,j+1}$, $\hat{\mathcal{A}}$ and $\hat{\mathcal{B}}_{i+1,j+1}$ in the following way:
\begin{equation}
T\hat \Omega_{i+1,j+1}T^{\top}= \Omega_{i+1,j+1},\ \ T\hat{\mathcal{A}}T^{\top}=\bar{\mathcal{A}}, \ \ T\hat{ \mathcal{B}}_{i+1,j+1} T^{\top}=\bar {\mathcal{B}}_{i+1,j+1},
\label{e21}
\end{equation}
where $\Omega_{i+1,j+1}$, $\bar{\mathcal{A}}$ and $\bar{ \mathcal{B}}_{i+1,j+1}$ were defined in (\ref{e9}).

Multiply the system ({\ref{es8}}) on the left by the matrix  $T$ and write the unknown vector $\bar w_{i+1,j+1}$ as $\bar w_{i+1,j+1}=T^{\top}z_{i+1,j+1}$. Then, taking into account (\ref{e20}) and (\ref{e21}), we obtain from (\ref{es8}) the following system:
\begin{equation}
\bar \Omega_{i+1,j+1}z_{i+1,j+1}=\tau g_{i+1,j+1}-(\bar{\mathcal{A}}+\tau\tilde \sigma_{3})z_{i+1,j}-r(\bar{\mathcal{B}}_{i+1,j+1}+\tau\tilde \sigma_{4})z_{i,j+1},
\label{e22}
\end{equation}
where
\[
\bar \Omega_{i+1,j+1}=\Omega_{i+1,j+1}+h\tilde\sigma_{1}+\tau\tilde\sigma_{2}, g_{i+1,j+1}=T\tilde f_{i+1,j+1}
\]
and $\tilde\sigma_{k}\equiv T\sigma_{k}T^{\top}$.

Next we suppose that in each node of the grid $U_{\Delta}$ inequalities (\ref{e11})  from the lemma \ref{l2} are valid and $r$  being a ratio of two steps is a constant. In virtue of the lemma \ref{l2}, the matrix $\Omega_{i+1,j+1}$ is nonsingular on the grid $U_{\Delta}$. Thus, the matrix $\bar\Omega_{i+1,j+1}$ in (\ref{es8}) can be presented in the following form:
\[
\bar\Omega_{i+1,j+1}=\Omega_{i+1,j+1}M_{i+1,j+1},
\]
where $M_{i+1,j+1}=E_{\tilde n}+\Omega^{-1}_{i+1,j+1}(h\tilde\sigma_{1}+\tau\tilde\sigma_{2})$. Since, by assumption, the elements of matrix coefficients of the system (\ref{e1}) belong to $C^{2}(U)$, then in virtue of the theorem 1  in the work \cite{verb}, the elements of  $J(x,t)$ also belong to $C^{2}(U)$.  This means that there is such a constant  $\mathcal{K}$, for which the condition  $\|\Omega^{-1}_{i+1,j+1}\|\leq \mathcal{K}$ is fulfilled.

In virtue of the theorem  from (\cite{lank}, p. 195) there exists the inverse matrix ${M^{-1}_{i+1,j+1}}$. It can be presented in the form
$M^{-1}_{i+1,j+1}=E_{\tilde n}+\tilde M_{i+1,j+1}h$, where
\[
\tilde M_{i+1,j+1}=\sum\limits_{\nu=1}^{\infty}(-1)^{\nu}\left [\Omega^{-1}_{i+1,j+1}(\sigma_{1}+r\tilde\sigma_{2})\right ]^{\nu}h^{\nu-1}
\]
Therefore we see that the matrix $\bar\Omega^{-1}_{i+1,j+1}$ exists. It is worth remarking, that therefore we proved the solvability of the difference scheme (\ref{e7}) in each node of the grid $U_{\Delta}$ for sufficiently small steps $h$ and $\tau$.

To proceed a transformation of the difference scheme (\ref{e7}), we  multiply (\ref{e22}) on the left by the matrix  $\bar\Omega_{i+1,j+1}^{-1}$ to get
\begin{equation}
z_{i+1,j+1}=\tilde g_{i+1,j+1}-(F_{i+1,j+1}+\bar \epsilon_{i+1,j+1}(h,\tau))z_{i+1,j}-r(K_{i+1,j+1}+\hat\epsilon_{i+1,j+1}(h,\tau)) z_{i,j+1},
\label{e23}
\end{equation}
where $\tilde g_{i+1,j+1}$ is the $\tilde n$-dimensional vector given by  $\tilde g_{i+1,j+1}=\tau M_{i+1,j+1}^{-1}\Omega_{i+1,j+1}^{-1}g_{i+1,j+1}$;
\[
F_{i+1,j+1}=\Omega_{i+1,j+1}^{-1}\bar{\mathcal{A}};\;\;\; K_{i+1,j+1}=\Omega_{i+1,j+1}^{-1}\bar{\mathcal{B}}_{i+1,j+1};
\]
and
\[
\bar\epsilon_{i+1,j+1}(h,\tau)=h\tilde M_{i+1,j+1}\Omega_{i+1,j+1}^{-1}\bar{\mathcal{A}}+\tau \Omega_{i+1,j+1}^{-1}\tilde \sigma_{3}+h\tau\tilde M_{i+1,j+1}\Omega_{i+1,j+1}^{-1}\tilde\sigma_{3},
\]
\[
\hat\epsilon_{i+1,j+1}(h,\tau)=h\tilde M_{i+1,j+1}\Omega_{i+1,j+1}^{-1}\bar{\mathcal{B}}_{i+1,j+1}+\tau \Omega_{i+1,j+1}^{-1}\tilde\sigma_{4}+h\tau\tilde M_{i+1,j+1}\Omega_{i+1,j+1}^{-1}\tilde\sigma_{4}
\]
being  $\tilde n\times\tilde n$ matrices.

Taking into account (\ref{e9}), the system (\ref{e23}) is decomposed  into three subsystems
\begin{eqnarray}
\bar w_{i+1,j+1}^{k}&=&\tilde g_{i+1,j+1}^{k}-F_{i+1,j+1}^{k}\bar w_{i+1,j}^{k}-rK_{i+1,j+1}^{k}\bar w_{i,j+1}^{k} \nonumber \\
&&-\sum_{l=1}^{3}\bar\epsilon_{i+1,j+1}^{k,l}(h,\tau)\bar w_{i+1,j}^{l}-\sum_{l=1}^{3}\hat\epsilon_{i+1,j+1}^{k,l}(h,\tau) \bar w_{i,j+1}^{l}.
\label{e24}
\end{eqnarray}
Here
\[
F_{i+1,j+1}^{k}=\left [\Omega_{i+1,j+1}^{k}\right ]^{-1}\mathcal{A}^{k},\;\;\;
K_{i+1,j+1}^{k}=\left [\Omega_{i+1,j+1}^{k}\right ]^{-1}\mathcal{B}_{i+1,j+1}^{k},
\]
while $\bar\epsilon_{i+1,j+1}^{k,l}(h,\tau)$ and $\hat\epsilon_{i+1,j+1}^{k,l}(h,\tau)$ are matrix blocks of 
\[
\bar\epsilon_{i+1,j+1}(h,\tau)=\left(\bar\epsilon_{i+1,j+1}^{k,l}(h,\tau)\right)_{1}^{3}\;\;\;\mbox{and}\;\;\; \hat\epsilon_{i+1,j+1}(h,\tau)=\left(\hat\epsilon_{i+1,j+1}^{k,l}(h,\tau)\right)_{1}^{3},
\] 
respectively.  The sizes of these blocks correspond to the decomposition of $z_{i+1, j+1}$. Finally, free terms $\tilde g_{i+1,j+1}^{k}$  are blocks of the vector
\[
\tilde g_{i+1,j+1}=(\tilde g_{i+1,j+1}^{1},\;\;\; \tilde g_{i+1,j+1}^{2},\;\;\; \tilde g_{i+1,j+1}^{3})^{\top}
\]
of the sizes $m_{1}m_{2}d$, $m_{1}m_{2}l$ and $m_{1}m_{2}p$, respectively, and the vectors $\bar w_{i+1,j}^{k}$ and $\bar w_{i,j+1}^{k}$ have the following form
\[
\bar w_{i+1,j}^{k}=(e_{m_{2}}\otimes w_{i+1,j}^{k},e_{m_{2}}\otimes w_{i+2,j}^{k}, \dots, e_{m_{2}}\otimes w_{i+m_{1},j}^{k})^{\top},
\]
\[
\bar w_{i,j+1}^{k}=e_{m_{1}}\otimes (w_{i,j+1}^{k},w_{i,j+2}^{k}, \dots, w_{i,j+m_{2}}^{k})^{\top}.
\]

Using (\ref{el19}), we transform the product $F_{i+1,j+1}^{1}\bar w_{i+1,j}^{1}$ in  the first equation of the system (\ref{e24}) to get
\begin{eqnarray}
F_{i+1,j+1}^{1}\bar w_{i+1,j}^{1}&=&\left(E_{m_{1}}\otimes\gamma_{m_{2}}\otimes E_{d}+r \bar\gamma_{m_{1}}\otimes E_{m_{2}}\otimes {J}_{i+1,j+1}\right )^{-1}\left (E_{m_{1}}\otimes \gamma^{0}_{m_{2}}\otimes E_{d}\right )\bar w_{i+1,j}^{1} \nonumber \\
&=&\left (E_{m_{1}m_{2}d}+r \bar\gamma_{m_{1}}\otimes \gamma_{m_{2}}^{-1}\otimes {J}_{i+1,j+1}\right )^{-1}\bar w_{i+1,j}^{1}+O(\tau^{m_{2}}).
\nonumber
\end{eqnarray}
Let $R$ be a  $m_{1}\times m_1$ matrix transforming $\bar\gamma_{m_{1}}$ to the normal Jordan form. Remark that the matrix $\bar\gamma_{m_{1}}$ is simple for any value of $m_{1}$, that is,
\[
R\bar\gamma_{m_{1}}R^{-1}=\bar\gamma_{m_{1}}^{*},\;\;\; \mbox{where}\;\;\;  \bar\gamma_{m_{1}}^{*}={\rm diag}\left\{\xi_{\bar\gamma_{m_{1}}}^{1},\xi_{\bar\gamma_{m_{1}}}^{2},\dots,\xi_{\bar\gamma_{m_{1}}}^{m_{1}}\right\}.
\]
Then
\begin{equation}
F_{i+1,j+1}^{1}\bar w_{i+1,j}^{1}=\left(R^{-1}\otimes E_{m_{2}d}\right)\tilde R\left(R\otimes E_{m_{2}d}\right)\bar w_{i+1,j}^{1}+O(\tau^{m_{2}}),
\label{e25}
\end{equation}
with
$\tilde R={\rm diag}\{\tilde R_{11},\tilde R_{22},\dots,\tilde R_{m_{1}m_{1}} \}$,
where
\[
\tilde R_{ss}=\left(E_{m_{2}d}+r\xi_{\bar\gamma_{m_{1}}}^{s}\left(\gamma_{m_{2}}^{-1}\otimes J_{i+1,j+1}\right) \right )^{-1}
\]
for $s=1,\ldots, m_1$. In virtue of the lemma \ref{l3} we have
\begin{equation}
\tilde R(R\otimes E_{m_{2}d})\bar w_{i+1,j}^{1}=\mathcal{D}^{1}_{i+1,j+1}(R\otimes E_{m_{2}d})\bar w_{i+1,j}^{1},
\label{e26}
\end{equation}
where
\[
\mathcal{D}^{1}_{i+1,j+1}={\rm diag} \left\{\exp\left(-r\xi_{\bar\gamma_{m_{1}}}^{1}J_{i+1,j+1}\right),\exp\left(-2r\xi_{\bar\gamma_{m_{1}}}^{1}J_{i+1,j+1}\right),\dots,\exp\left(-m_{2}r\xi_{\bar\gamma_{m_{1}}}^{1}J_{i+1,j+1}\right), \right.
\]
\[
\exp\left(-r\xi_{\bar\gamma_{m_{1}}}^{2}J_{i+1,j+1}\right),\exp\left(-2r\xi_{\bar\gamma_{m_{1}}}^{2}J_{i+1,j+1}\right),\dots,\exp\left(-m_{2}r\xi_{\bar\gamma_{m_{1}}}^{2}J_{i+1,j+1}\right),\dots
\]
\[
\left.
\dots
,\exp\left(-r\xi_{\bar\gamma_{m_{1}}}^{m_{1}}J_{i+1,j+1}\right),\exp\left(-2r\xi_{\bar\gamma_{m_{1}}}^{m_{1}}J_{i+1,j+1}\right),\dots,\exp\left(-m_{2}r\xi_{\bar\gamma_{m_{1}}}^{m_{1}}J_{i+1,j+1}\right)\right\}.
\]
Then, taking into account (\ref{e26}), the relation (\ref{e25}) takes the form
\begin{equation}
F_{i+1,j+1}^{1}\bar w_{i+1,j}^{1}=(R^{-1}\otimes E_{m_{2}d})\mathcal{D}^{1}_{i+1,j+1}(R\otimes E_{m_{2}d})\bar w_{i+1,j}^{1}+O(\tau^{m_{2}}).
\label{e28}
\end{equation}
Next we consider block  permutation matrix $\mathcal{T}=(T_{i,j})$, where $i=1,\ldots, m_2$ and $j=1,\ldots, m_1$, where each block  $T_{i,j}$ itself consists of blocks $T_{i,j}=({\tilde T}_{k,s})$, for $k=1,\ldots, m_1$ and $s=1,\ldots, m_2$. The blocks ${\tilde T}_{k,s}$ are defined by
\[
{\tilde T}_{k,s}=\left \{\begin{array}{l}
E_{d},\ \mbox{for}\ k=j,\ s=i,\\
\mathcal{O}_{d},\  \mbox{for}\ k\neq j \ \mbox{or} \  s\neq i
\end{array}\right . .
\]
Remark, that the matrix $\mathcal{T}$ is in fact a square orthogonal matrix of the order $m_{1}m_{2}d$. Also it is worth remarking, that
\begin{equation}
\mathcal{T}\left(R^{-1}\otimes E_{m_{2}d}\right)\mathcal{D}^{1}_{i+1,j+1}\left(R\otimes E_{m_{2}d}\right)\mathcal{T}^{\top}=\tilde{\mathcal {D}}^{1}_{i+1,j+1},
\label{e29}
\end{equation}
where
\[
\tilde{\mathcal{D}}^{1}_{i+1,j+1}={\rm diag}\left\{\exp\left(-r\bar\gamma_{m_{1}}\otimes J_{i+1,j+1}\right),\exp\left(-2r\bar\gamma_{m_{1}}\otimes J_{i+1,j+1}\right),\dots,\exp\left(-m_{2}r\bar\gamma_{m_{1}}\otimes J_{i+1,j+1}\right)\right\}.
\]
Then the relation (\ref{e28}) taking into account (\ref{e29}) becomes
\begin{equation}
F_{i+1,j+1}^{1}\bar w_{i+1,j}^{1}=\mathcal{T}^{\top}\tilde{\mathcal{D}}^{1}_{i+1,j+1}\mathcal{T}\bar w_{i+1,j}^{1}+O(\tau^{m_{2}}).
\label{e30}
\end{equation}
Therefore we finished a transformation of $F_{i+1,j+1}^{1}\bar w_{i+1,j}^{1}$ and now let us turn to $rK_{i+1,j+1}^{1}\bar w_{i,j+1}^{1}$ in  the first equation  of system (\ref{e24}). Using (\ref{e9}) and taking into account (\ref{el19}) we obtain
\begin{eqnarray}
rK_{i+1,j+1}^{1}\bar w_{i,j+1}^{1}&=&r\left (E_{m_{1}}\otimes\gamma_{m_{2}}\otimes E_{d}+r \bar\gamma_{m_{1}}\otimes E_{m_{2}}\otimes {J}_{i+1,j+1}\right )^{-1}\left (\gamma_{m_{1}}^{0}\otimes E_{m_{2}}\otimes J_{i+1,j+1}\right )\bar w_{i,j+1}^{1} \nonumber \\
&=&\left (E_{m_{1}m_{2}d}+\frac{1}{r} \bar\gamma_{m_{1}}^{-1}\otimes \gamma_{m_{2}}\otimes {J}_{i+1,j+1}^{-1}\right )^{-1}\bar w_{i,j+1}^{1}+O(h^{m_{1}}). \nonumber
\end{eqnarray}
Again, let $R_{1}$ be  $m_{2}\times m_{2}$ matrix, transforming $\gamma_{m_{2}}$ to normal Jordan form $\gamma_{m_{2}}^{*}$, that is,
\[
R_{1}\gamma_{m_{2}}R_{1}^{-1}=\gamma_{m_{2}}^{*},\;\;\; \mbox{where}\;\;\;  \gamma_{m_{2}}^{*}={\rm diag}\left\{\xi_{\gamma_{m_{2}}}^{1},\xi_{\gamma_{m_{2}}}^{2},\dots,\xi_{\gamma_{m_{2}}}^{m_{2}}\right\}.
\]
Then
\begin{equation}
rK_{i+1,j+1}^{1}\bar w_{i,j+1}^{1}=\left(E_{m_{1}}\otimes R_{1}^{-1}\otimes E_{d}\right)\tilde R_{1}\left(E_{m_{1}}\otimes R_{1}\otimes E_{d}\right)\bar w_{i,j+1}^{1}+O(h^{m_{1}}),
\label{e31}
\end{equation}
where
\[
\tilde R_{1}=\left ( E_{m_{1}m_{2}d}+\frac{1}{r}\bar\gamma_{m_{1}}^{-1}\otimes\gamma_{m_{2}}^{*}\otimes J^{-1}_{i+1,j+1} \right )^{-1}.
\]
Remark, that
\begin{equation}
\mathcal{T}\tilde R_{1}\mathcal{T}^{\top}=\bar R,
\label{e32}
\end{equation}
with $\bar R={\rm diag}\{\bar R_{11}, \bar R_{22},\dots, \bar R_{m_{2}m_{2}}\}$, where
\[
\bar R_{ss}\equiv\left( E_{m_{1}d}+\frac{1}{r}\xi_{m_{2}}^{s}\left(\gamma_{m_{1}}^{-1}\otimes J^{-1}_{i+1,j+1}\right)\right)^{-1},
\]
for $s=1,\ldots,m_2$. From (\ref{e31}) and (\ref{e32}),  taking into account lemma {\ref{l3}}, we get
\begin{equation}
rK_{i+1,j+1}^{1}\bar w_{i,j+1}^{1}=\left(E_{m_{1}}\otimes R_{1}^{-1}\otimes E_{d}\right)\mathcal{T}^{\top}\mathcal{D}^{2}_{i+1,j+1}\mathcal{ T}\left(E_{m_{1}}\otimes R_{1}\otimes E_{d}\right)\bar w_{i,j+1}^{1}+O(h^{m_{1}}),
\label{e33}
\end{equation}
where
\[
\mathcal{D}^{2}_{i+1,j+1}={\rm diag} \left\{ \exp\left( -\frac{\xi_{\gamma_{m_{2}}}^{1}}{r}J_{i+1,j+1}^{-1}\right),\exp\left( -\frac{2\xi_{\gamma_{m_{2}}}^{1}}{r}J_{i+1,j+1}^{-1}\right),\dots, \exp\left( -\frac{m_{1}\xi_{\gamma_{m_{2}}}^{1}}{r}J_{i+1,j+1}^{-1}\right),\right.
\]
\[
\exp\left( -\frac{\xi_{\gamma_{m_{2}}}^{2}}{r}J_{i+1,j+1}^{-1}\right),\exp\left( -\frac{\xi_{\gamma_{m_{2}}}^{2}}{r}J_{i+1,j+1}^{-1}\right),\dots, \exp\left( -\frac{m_{1}\xi_{\gamma_{m_{2}}}^{2}}{r}J_{i+1,j+1}^{-1}\right)
,\dots
\]
\[
\left.
\dots,
\exp\left( -\frac{\xi_{\gamma_{m_{2}}}^{m_{1}}}{r}J_{i+1,j+1}^{-1}\right),\exp\left( -\frac{2\xi_{\gamma_{m_{2}}}^{m_{1}}}{r}J_{i+1,j+1}^{-1}\right),\dots, \exp\left( -\frac{m_{1}\xi_{\gamma_{m_{2}}}^{m_{1}}}{r}J_{i+1,j+1}^{-1}\right)\right\}.
\]
Remark, that
\begin{equation}
\left(E_{m_{1}}\otimes R_{1}^{-1}\otimes E_{d}\right)\mathcal{T}^{\top}\mathcal{D}^{2}_{i+1,j+1}\mathcal{T}\left(E_{m_{1}}\otimes R_{1}\otimes E_{d}\right)=\tilde{\mathcal{D}}^{2}_{i+1,j+1},
\label{e34}
\end{equation}
where
\[
\tilde{\mathcal{D}}^{2}_{i+1,j+1}={\rm diag}\left\{\exp\left(-\frac{1}{r}\gamma_{m_{2}}\otimes J_{i+1,j+1}^{-1}\right),\exp\left(-\frac{2}{r}\gamma_{m_{2}}\otimes J_{i+1,j+1}^{-1}\right),\dots,\exp\left(-\frac{m_{1}}{r}\gamma_{m_{2}}\otimes J_{i+1,j+1}^{-1}\right)\right\}.
\]
Therefore, from (\ref{e33}) and (\ref{e34}) we get
\begin{equation}
rK_{i+1,j+1}^{1}\bar w_{i,j+1}^{1}=\tilde{\mathcal{D}}^{2}_{i+1,j+1}\bar w_{i,j+1}^{1}+O(h^{m_{1}}).
\label{e35}
\end{equation}
To transform the coefficients of the second and third equation of the system (\ref{e24}), it is enough  to make use of (\ref{e9}) and (\ref{el19}). In a result we obtain
\[
F_{i+1,j+1}^{2}\bar w_{i+1,j}^{2}=0,\ \
rK_{i+1,j+1}^{2}\bar w_{i,j+1}^{2}=\bar w_{i,j+1}^{2}+O(h^{m_{1}}),
\]
\begin{equation}
F_{i+1,j+1}^{3}\bar w_{i+1,j}^{3}=\bar w_{i+1,j}^{3}+O(\tau^{m_{2}}),\ \
rK_{i+1,j+1}^{3}\bar w_{i,j+1}^{3}=0.
\label{e36}
\end{equation}
Therefore, a transformation of all coefficients  of the system (\ref{e24}) is finished. Let us write down transformed difference scheme using (\ref{e30}), (\ref{e35}) and (\ref{e36}). We have
\begin{eqnarray}
\bar w_{i+1,j+1}^{1}&=&\tilde g_{i+1,j+1}^{1}-\mathcal{T}^{\top}\tilde{\mathcal{D}}^{1}_{i+1,j+1}\mathcal{T}\bar w_{i+1,j}^{1}-\tilde{\mathcal {D}}^{2}_{i+1,j+1}\bar w_{i,j+1}^{1}-
\sum_{l=1}^{3}\bar\epsilon_{i+1,j+1}^{1,l}(h,\tau)\bar w_{i+1,j}^{l} \nonumber \\
&&-\sum_{l=1}^{3}\hat\epsilon_{i+1,j+1}^{1,l}(h,\tau) \bar w_{i,j+1}^{l}+O(h^{m_{1}})+O(\tau^{m_{2}}), \label{e37}
\end{eqnarray}
\[
\bar w_{i+1,j+1}^{2}=\tilde g_{i+1,j+1}^{2}-\bar w_{i,j+1}^{2}
-\sum_{l=1}^{3}\bar\epsilon_{i+1,j+1}^{2,l}(h,\tau)\bar w_{i+1,j}^{l}-\sum_{l=1}^{3}\hat\epsilon_{i+1,j+1}^{2,l}(h,\tau) \bar w_{i,j+1}^{l}+O(h^{m_{1}}),
\]
\[
\bar w_{i+1,j+1}^{3}=\tilde g_{i+1,j+1}^{3}-\bar w_{i+1,j}^{3}-
\sum_{l=1}^{3}\bar\epsilon_{i+1,j+1}^{3,l}(h,\tau)\bar w_{i+1,j}^{l}-\sum_{l=1}^{3}\hat\epsilon_{i+1,j+1}^{3,l}(h,\tau) \bar w_{i,j+1}^{l}+O(\tau^{m_{2}}).
\]
The system (\ref{e37}) represents a canonical form of the difference scheme (\ref{e9}). In this form one immediately sees spectral characteristics of matrix coefficients. This information is necessary to prove a uniform boundedness of the grid solution in the domain  $U$.
Therefore the goal of this section is achieved and we are in position to turn to the main part of the paper.

\section{The proof of a stability of the difference scheme}
\label{sect6}

In this section we prove the following.
\begin{theorem}
Let in the difference scheme (\ref{e7}) eigenvalues of the matrices $\bar\gamma_{m_{1}}$, $\gamma_{m_{2}}$ and $J_{i+1,j+1}$ with given $m_{1}$, $m_{2}$ satisfy  the condition (\ref{e11}) in the domain $U_{\Delta}$ and $\xi_{{J}_{i+1,j+1}}^{s}>0$ $\forall\ i,j$ and  $s=1,\ldots, k$. Let $r$ being a ratio of two steps is a constant.  Then the difference scheme (\ref{e7}) is absolutely stable with respect to initial-boundary data and right-hand side. Following estimation for its solution
\begin{equation}
\|v_{i+1,j+1}\|_{C(U_{\Delta})}\leq \mathcal{M}_{1}\|f_{i+1,j+1}\|_{C(U_{\Delta})}+\mathcal{M}_{2}\|\phi_{i+1}\|_{C(U_{\Delta})}+\mathcal{ M}_{3}\|\psi_{j+1}\|_{C(U_{\Delta})}
\label{et1}
\end{equation}
is valid, where $\mathcal{M}_k$ are constants and $i=1,\ldots, n_{1}-1$, $j=1,\ldots, n_{2}-1$.
\label{t2}
\end{theorem}
{\sc Proof.}
Remark that the existence of a unique solution of the difference scheme (\ref{e7}) was in fact  proven  in previous section together with its transformation. It remains to prove the second part of the theorem.

The difference scheme (\ref{e7}) in the previous section was transformed to the form (\ref{e37}). To describe the structure of difference scheme (\ref{e37}), we rewrite it in matrix form
\begin{equation}
(L+\mathcal{L}(h,\tau))\mathcal{V}= G,
\label{e38}
\end{equation}
where $\mathcal{V}$ is unknown $n_{1}n_{2}\tilde n$-dimensional vector, consisting of three block components $\mathcal{V}=(\mathcal{V}^{1}\ \mathcal{V}^{2}\ \mathcal{V}^{3})$, where
\[
\mathcal{V}^{k}=\left(\bar w^{k}_{1,1},\ \bar w^{k}_{1,2},\dots,\bar w^{k}_{1,n_{2}},\bar w^{k}_{2,1},\ \bar w^{k}_{2,2},\dots,\bar w^{k}_{2,n_{2}},\dots, \bar w^{k}_{n_{1},1},\ \bar w^{k}_{n_{1},2},\dots,\bar w^{k}_{n_{1},n_{2}}\right)^{\top}.
\]

Let us describe in detail system (\ref{e38}). In (\ref{e38}),  $L$ is a square block diagonal matrix of the order $n_{1}n_{2}\tilde n$, that is,  $L={\rm diag}\{L^{1},L^{2},L^{3}\}$. In turn each block $L^{k}$ is block two-diagonal matrix $L^{k}=(L_{i,j}^{k})$ for $i=1,\ldots, n_{1}$ and $j=1,\ldots, n_{2}$. The blocks $L_{i,j}^{k}$ at the main diagonal, also have the block two-diagonal form
\begin{equation}
L_{i,j}^{k}=\left ( \begin{array}{cccccc}
E_{s} & \mathcal{O}_{s} &  \mathcal{O}_{s} & \dots & \mathcal{O}_{s} & \mathcal{O}_{s}\\
\mathcal{F}_{i,2}^{k} & E_{s} &  \mathcal{O}_{s} & \dots & \mathcal{O}_{s} & \mathcal{O}_{s}\\
\mathcal{O}_{s} & \mathcal{F}_{i,3}^{k} & E_{s} & \dots & \mathcal{O}_{s} & \mathcal{O}_{s}\\
\vdots & \vdots &  \vdots & \ddots & \vdots & \vdots\\
\mathcal{O}_{s} & \mathcal{O}_{s} &  \mathcal{O}_{s} & \dots & E_{s} & \mathcal{O}_{s}\\
\mathcal{O}_{s} & \mathcal{O}_{s} &  \mathcal{O}_{s} & \dots & \mathcal{F}_{i,n_{2}}^{k} & E_{s}\\
\end{array}
\right ).
\label{e39}
\end{equation}
Here the parameter $s$ take values $d,\ l$ and $p$ for $k=1,\ 2$ and $3$, respectively. The blocks $\mathcal{F}_{i,j}^{k}$  of the matrix (\ref{e39}) have the form
\begin{equation}
\mathcal{F}^{1}_{i,j}=\mathcal{T}^{\top}\tilde{\mathcal{D}}^{1}_{i,j}\mathcal{T},\ \ \  \mathcal{F}^{2}_{i,j}=\mathcal{O}_{l},\ \ \ \mathcal{F}^{3}_{i,j}=E_{p}.
\label{e40}
\end{equation}
The blocks $L_{i,j}^{k}$, situated under the main diagonal, that is,  for $i=2,\ldots, n_{1}$ and $j=i-1$, have the following block diagonal form:
\begin{equation}
L_{i,j}^{k}={\rm diag}\left\{\mathcal{K}_{i,1}^{k},\ \mathcal{K}_{i,2}^{k}, \dots , \mathcal{K}_{i,n_{2}}^{k}\right\},\ \ \
\mathcal{K}_{i,j}^{1}=\tilde{\mathcal{D}}^{2}_{i,j},\ \ \ \mathcal{K}_{i,j}^{2}=E_{l},\ \ \ \mathcal{K}_{i,j}^{3}=\mathcal{O}_{p}.
\label{e41}
\end{equation}
All the rest blocks $L_{i,j}^{k}$ are zero ones of corresponding sizes.

The vector $G$ in the system (\ref{e38}) is as follows:
\begin{equation}
G=\bar g+(S+\mathcal{S}(h,\tau))w_{0}+(Q+\mathcal{Q}(h,\tau))w^{0}+O(h^{m_{1}})+O(\tau^{m_{2}}),
\label{e42}
\end{equation}
with $\bar g$ are $n_{1}n_{2}\tilde n$-dimensional vector, consisting of blocks $\bar g=(g^{1},\ g^{2},\ g^{3})^{\top}$ of the sizes $n_{1}n_{2}d$, $n_{1}n_{2}l$ and $n_{1}n_{2}p$, respectively,
\[
g^{k}=(\tilde g^{k}_{1,1},\ \tilde g^{k}_{1,2},\dots, \tilde g^{k}_{1,n_{2}},\tilde g^{k}_{2,1},\ \tilde g^{k}_{2,2},\dots,\tilde g^{k}_{2,n_{2}},\dots, \tilde g^{k}_{n_{1},1},\ \tilde g^{k}_{n_{1},2},\dots,\tilde g^{k}_{n_{1},n_{2}})^{\top};
\]
$w_{0}$ and $w^{0}$ are known vectors of the size $n_{1}n_{2}\tilde n$ the elements of which are defined by initial-boundary data (\ref{e2}), namely, $w_{0}=(w_{0}^{1},\ w_{0}^{2},\ w_{0}^{3})^{\top}$,  where
\[
w_{0}^{k}=\left(\bar w_{1,0}^{k},\ \bar w_{2,0}^{k},\dots, \bar w_{n_{1},0}^{k}\right)^{\top}\otimes e_{n_{2}}
\]
and
$w^{0}=\left(w^{0,1},\ w^{0,2},\ w^{0,3}\right)^{\top}$, where
\[
w^{0,k}=e_{n_{2}}\otimes \left(\bar w_{0,1}^{k},\ \bar w_{0,2}^{k},\dots, \bar w_{0,n_{2}}^{k}\right)^{\top}.
\]
The matrices $S$ and $Q$ are known square block diagonal ones of the order $n_{1}n_{2}\tilde n$, that is,
\[
S={\rm diag}\left\{S^{1},\ S^{2},\ S^{3}\right\}\;\;\;\mbox{and}\;\;\; Q={\rm diag}\left\{Q^{1},\ Q^{2},\ Q^{3}\right\}.
\]
Each square block $S^{k}$ and $Q^{k}$ has the order $n_{1}n_{2}d$, $n_{1}n_{2}l$ and $n_{1}n_{2}p$, corresponding to the value of  $k$.
Each block $S^{k}$ is also block diagonal matrix
\[
S^{k}={\rm diag}\left\{S_{1,1}^{k},\ S_{2,2}^{k},\dots,S_{n_{1},n_{1}}^{k}\right\},
\]
where $S_{i,i}^{k}$ are square blocks of the orders $n_{2}d$, $n_{2}l$ and $n_{2}p$, respectively $k$, of the following form:
\[
S_{i,i}^{k}={\rm diag}\left\{-\mathcal{F}_{i,1}^{k},\ \mathcal{O}_{s},\  \mathcal{O}_{s},\dots, \mathcal{O}_{s}\right\},
\]
for $i=1,\ldots, n_{1}$,where $s$ takes values  $d$, $l$ and $p$, respectively  $k$.
Each block $Q^{k}$ is block diagonal matrix
\[
Q^{k}={\rm diag}\left\{Q_{1,1}^{k},\ \mathcal{O}_{s},\ \mathcal{O}_{s},\dots, \mathcal{O}_{s}\right\},
\]
where
\[
Q_{1, 1}^{k}={\rm diag}\left\{-\mathcal{K}_{1, 1}^{k},\ -\mathcal{K}_{1, 2}^{k},\dots, -\mathcal{K}_{1, n_{2}}^{k}\right\}.
\]

Finally, the matrices $\mathcal{L}(h,\tau)$, $\mathcal{S}(h,\tau)$ and $\mathcal{Q}(h,\tau)$  in (\ref{e38}) and (\ref{e42}) are quadratic matrices of the order $n_{1}n_{2}\tilde n$ built from the blocks  $\bar\epsilon_{i+1,j+1}^{s,l}(h,\tau)$ and $\hat\epsilon_{i+1,j+1}^{s,l}(h,\tau)$.

To prove the stability property of difference scheme (\ref{e7}), there is a need  to estimate the norm  of unknown vector $\mathcal{V}$ in the system  (\ref{e38}). For this aim we must calculate $L^{-1}$. Since $L$ is a block diagonal matrix, and each its diagonal block is a block two-diagonal matrix, then we can easily write down  the matrix $L^{-1}$ in its explicit form. Each its diagonal block component is of the form
\begin{equation}
(L^{k})^{-1}=\Lambda^{k}\Phi^{k},\;\;\;\mbox{for}\;\;\; k=1, 2, 3,
\label{e43}
\end{equation}
where $\Lambda^{k}$ are matrices of the order  $n_{1}n_{2}d$, $n_{1}n_{2}l$ and $n_{1}n_{2}p$, respectively value of $k$, which have the form
\begin{equation}
\Lambda^{k}=\left ( \begin{array}{ccccc}
E_{s} & \mathcal{O}_{s} &  \mathcal{O}_{s} & \dots & \mathcal{O}_{s}\\
\Lambda^{k}_{2} & E_{s} &  \mathcal{O}_{s} & \dots & \mathcal{O}_{s}\\
\Lambda^{k}_{3}\Lambda^{k}_{2} & \Lambda^{k}_{3} & E_{s} & \dots & \mathcal{O}_{s}\\
\vdots & \vdots &  \vdots & \ddots & \vdots\\
\prod_{s_{1}=2}^{n_{1}}\Lambda^{k}_{s_{1}} & \prod_{s_{1}=3}^{n_{1}}\Lambda^{k}_{s_{1}} &  \prod_{s_{1}=4}^{n_{1}}\Lambda^{k}_{s_{1}} & \dots & E_{s}\\
\end{array}
\right ),
\label{e44}
\end{equation}
where $s$ take values $n_{2}d,\ n_{2}l$ and $n_{2}p$, respectively value of $k$.
Every block $\Lambda^{k}_{i}$, for $i=2,\ldots, n_{1}$ is defined as $\Lambda^{k}_{i}=(\Phi_{i}^{k})^{-1}Q_{i,1}^{k}$, where
\begin{equation}
(\Phi_{i}^{k})^{-1}=\left ( \begin{array}{ccccc}
E_{s} & \mathcal{O}_{s} &  \mathcal{O}_{s} & \dots & \mathcal{O}_{s}\\
-\mathcal{F}^{k}_{i,2} & E_{s} &  \mathcal{O}_{s} & \dots & \mathcal{O}_{s}\\
\mathcal{F}^{k}_{i,3}\mathcal{F}^{k}_{i,2} & -\mathcal{F}^{k}_{i,3} & E_{s} & \dots & \mathcal{O}_{s}\\
\vdots & \vdots &  \vdots & \ddots & \vdots\\
(-1)^{n_{2}+1}\prod_{s_{1}=2}^{n_{2}}\mathcal{F}^{k}_{i,s_{1}} & (-1)^{n_{2}+2}\prod_{s_{1}=3}^{n_{2}}\mathcal{F}^{k}_{i,s_{1}} &  (-1)^{n_{2}+3}\prod_{s_{1}=4}^{n_{2}}\mathcal{F}^{k}_{i,s_{1}} & \dots & E_{s}\\
\end{array}
\right ),
\label{e45}
\end{equation}
where $s$ take values $d,\ l$ and $p$, respectively value of $k$.
The matrices $\Phi^{k}$  in (\ref{e43}) have the form
\[
\Phi^{k}={\rm diag}\{(\Phi_{1}^{k})^{-1},\ (\Phi_{2}^{k})^{-1},\dots,(\Phi_{n_{1}}^{k})^{-1} \}.
\]

Estimate the norms of the  matrices $(L^{k})^{-1}$. For $k=2$ and $k=3$,  it follows from (\ref{e40}), (\ref{e41}) and (\ref{e43})-(\ref{e45})  that
$\Lambda^{k}_{i}=\mathcal{O}_{s}$  and $(L^{k})^{-1}=E_{n_{1}n_{2}s}$, where $s=d,l$ for $k=2,3$, respectively. Thus,
\begin{equation}
\|(L^{k})^{-1}\|_{C(U_{\Delta})}=1\ \ \ \mbox {for}\ \ \ k=2,3.
\label{e46}
\end{equation}

Let us estimate the norm of the first block component  $(L^{1})^{-1}$. It follows from (\ref{e43}) and (\ref{e44}) that
\[
\|(L^{1})^{-1}\|_{C(U_{\Delta})}\leq \|\Lambda^{1}\|_{C(U_{\Delta})}\|\Phi^{1}\|_{C(U_{\Delta})},
\]
\begin{equation}
\|\Lambda^{1}\|_{C(U_{\Delta})}\leq  1+\|\Lambda^{1}_{i}\|_{C(U_{\Delta})}+\|\Lambda^{1}_{i}\|^{2}_{C(U_{\Delta})}+\cdots+\|\Lambda^{1}_{i}\|^{n_{1}-1}_{C(U_{\Delta})}.
\label{e47}
\end{equation}
In virtue of (\ref{e45})
\begin{equation}
\|\Lambda^{1}_{i}\|_{C(U_{\Delta})}\leq \|(\Phi_{i}^{1})^{-1}\|_{C(U_{\Delta})} \|\mathcal{K}_{i,j}^{1}\|_{C(U_{\Delta})}.
\label{e48}
\end{equation}
In turn for the first multiplier in the right-hand side of the inequality (\ref{e48}) we have
\begin{equation}
\|(\Phi_{i}^{1})^{-1}\|_{C(U_{\Delta})}\leq
1+\|\mathcal{F}^{1}_{i,n_{2}}\|_{C(U_{\Delta})}+\|\mathcal{F}^{1}_{i,n_{2}}\mathcal{F}^{1}_{i,n_{2}-1}\|_{C(U_{\Delta})}+\cdots+\|\mathcal{ F}^{1}_{i,n_{2}}\mathcal{F}^{1}_{i,n_{2}-1}\cdots\mathcal{F}^{1}_{i,2}\|_{C(U_{\Delta})}.
\label{e49}
\end{equation}
From (\ref{e29}) and (\ref{e40}) we find the radius $\mu$ of the spectrum of the matrix $\mathcal{F}^{1}_{i,j}$ to get
\begin{equation}
\mu=\max\left \{\left |\exp\left(-k_{1}r\xi_{\bar\gamma_{m_{1}}}^{k_{2}}\xi_{J_{i,j}}^{k_{3}}\right)\right |\right \}
\label{e50}
\end{equation}
for $i=1,\ldots,n_{1}$, $j=1,\ldots, n_{2}$, $k_{1}=1,\ldots, m_{2}$, $k_{2}=1,\ldots,m_{1}$ and $k_{3}=1,\ldots, d$.
We have a pair of inequalities: ${\rm Re}(\xi_{\bar\gamma_{m_{1}}}^{k_{2}})>0\ \forall\ k_{2}$, which can be verified with the help of Routh-Hurwitz criterion and $\xi_{J_{i,j}}^{k_{3}}>0\ \forall\ k_{3}$, which is valid by condition of our theorem. Then,  from (\ref{e50}) we get $\mu<1$. Applying theorem on spectral decomposition for a power function of matrix, namely,  $(\mathcal{F}^{1}_{i,j})^{n_{2}}$ given in (\cite{lank}, p. 155) or theorem from \cite{vorob}, we get
\begin{equation}
\|(\mathcal{F}^{1}_{i,j})^{n_{2}}\|_{C(U_{\Delta})}\to 0,\ \ \mbox{at}\ \ n_{2}\to\infty.
\label{e51}
\end{equation}
From (\ref{e51}) it follows, that there exists such a value of the power  $\tilde m$, for which
\begin{equation}
\|(\mathcal{F}^{1}_{i,j})^{\tilde m}\|_{C(U_{\Delta})}<1,
\label{e52}
\end{equation}
where $m$ is a suitable constant. Let us consider the product $\prod_{j=2}^{\tilde m+1}\mathcal{F}^{1}_{i,j}$. With Taylor's formula we obtain
\begin{equation}
\prod_{j=2}^{\tilde m+1}\mathcal{F}^{1}_{i,j}=(\mathcal{F}^{1}_{i,2})^{\tilde m}+O(\tau).
\label{e53}
\end{equation}
Denote $\chi=\left|\left|\prod_{j=2}^{\tilde m+1}\mathcal{F}^{1}_{i,j}\right|\right|_{C(U_{\Delta})}$. From (\ref{e52}) and (\ref{e53}) it follows that $\chi<1$.
Hence, from (\ref{e49}) we get
\begin{equation}
\|(\Phi_{i}^{1})^{-1}\|_{C(U_{\Delta})}\leq{\eta}/{(1-\chi)},
\label{e55}
\end{equation}
where
\[
\eta=1+\|\mathcal{F}^{1}_{i,j}\|_{C(U_{\Delta})}+\|\mathcal{F}^{1}_{i,j}\|^{2}_{C(U_{\Delta})}+\cdots+\|\mathcal{F}^{1}_{i,j}\|^{\tilde m-1}_{C(U_{\Delta})}.
\]

Let us estimate a norm of the matrix $\mathcal{K}^{1}_{i,j}$. We learn from \cite{bulg} that for a constant $n\times n$ matrix $\tilde A$ spectrum of which lies strictly in the left halfplane
 it is valid following inequality
\begin{equation}
\|\exp(t\tilde A)\|_{2}\leq \sqrt{c}\exp(-\kappa t),\ \ t\geq 0,
\label{e56}
\end{equation}
where $\kappa=1/(2\|X\|_{2})$, $c=\|X^{-1}\|_{2}\|X\|_{2}$ and  $X$ is Hermitian $n\times n$ matrix, being a solution of Lyapunov equation
$XA+A^{*}X=-E_{n}$ and $\|\cdot\|_{2}$ ia a norm of a matrix, agreed with  Hermitian norm. The constants $c$ and $\kappa$  are clarified in \cite{nechep}.

From the above mentioned assumptions with respect to eigenvalues of the matrix $J(x,t)$ and from (\ref{e34}),(\ref{e41}) and (\ref{e56}) it follows, that for the matrix $\mathcal{K}^{1}_{i,j}$ it is valid the following estimation:
\begin{equation}
\|\mathcal{K}^{1}_{i,j}\|_{C(U_{\Delta})}\leq \sqrt{\tilde c}\exp(-\tilde \kappa /r),
\label{e57}
\end{equation}
where $\tilde c$ and $\tilde \kappa$ are constants. From (\ref{e57}) it follows that, reducing  $r$, one can achieve a sufficient smallness of the elements of  $\mathcal{K}^{1}_{i,j}$ including that the inequality
\begin{equation}
\frac{\eta}{1-\chi}\|\mathcal{K}^{1}_{i,j}\|_{C(U_{\Delta})}<1.
\label{e58}
\end{equation}
will be valid.
So, from (\ref{e48}), (\ref{e55}) and (\ref{e58}) it follows that there are such values of the steps $\tau$ and $h$ for which the inequality $\|\Lambda^{1}_{i}\|_{C(U_{\Delta})}<1$ holds. Taking into account (\ref{e47}), we get
\begin{equation}
\|\Lambda^{1}\|_{C(U_{\Delta})}<\hat c,\ \ \ \mbox{where}\ \ \ \hat{c}=1/(1-\|\Lambda_{i}^{1}\|_{C(U_{\Delta})}).
\label{e59}
\end{equation}
Therefore, from (\ref{e47}) and (\ref{e55}) it follows that
\begin{equation}
\|(L^{1})^{-1}\|_{C(U_{\Delta})}<\bar c,\ \ \ \mbox{where}\ \ \ \bar c={\hat c\eta}/{(1-\chi)}.
\label{e60}
\end{equation}
These inequalities (\ref{e46}) and  (\ref{e60}) mean uniform boundedness of the matrix  $L^{-1}$ in (\ref{e38}).
From (\ref{e26}) and (\ref{e40}) it follows that to estimate a norm of the matrix $\mathcal{F}_{i,j}^{1}$, we can use the inequality (\ref{e56}). In a result, we get
\begin{equation}
\|\mathcal{F}_{i,j}^{1}\|_{C(U_{\Delta})}\leq \sqrt{\epsilon}\exp(-\bar \kappa r),
\label{e61}
\end{equation}
where $\epsilon$ and $\bar \kappa$ are some constants. Then the matrices $S$ and $Q$ in (\ref{e42}), in virtue of (\ref{e57}) and (\ref{e61}), are also are  uniformly boundness at the grid $U_{\Delta}$. That is, for sufficiently small values of $h$ and $\tau$ there are such constants  $\rho_{1}$ and $\rho_{2}$, for which
\begin{equation}
\|S+\mathcal{S}(h,\tau)\|_{C(U_{\Delta})}\leq\rho_{1}, \ \ \ \|Q+\mathcal{Q}(h,\tau)\|_{C(U_{\Delta})}\leq\rho_{2}.
\label{e62}
\end{equation}
From (\ref{e38}) and (\ref{e42}) we obtain
\begin{equation}
\mathcal{V}= (E_{n_{1}n_{2}\tilde n}+L^{-1}\mathcal{L}(h,\tau))^{-1} L^{-1} \left (\bar g+(S+\mathcal{S}(h,\tau))w_{0}+(Q+\mathcal{ Q}(h,\tau))w^{0}+O(h^{m_{1}})+O(\tau^{m_{2}})\right ).
\label{e63}
\end{equation}
In virtue of a boundness of the matrix $L^{-1}$,  for sufficiently small $\tau$ and $h$, we have
\[
\|(E_{n_{1}n_{2}\tilde n}+L^{-1}\mathcal{L}(h,\tau))^{-1}\|_{C(U_{\Delta})}\leq \bar \eta,
\]
where $\bar \eta=1/(1-\tilde\eta)$ and   $\tilde\eta=\|L^{-1}\mathcal{L}(h,\tau)\|_{C(U_{\Delta})}$. Denote
\[
\mathcal{M}_{1}=\tau\rho\bar\eta\mathcal{K}\|\tilde P\|_{C(U_{\Delta})}/(1-\|\tilde M_{i+1,j+1}\|_{C(U_{\Delta})}),
\]
\[
\mathcal{M}_{2}=\rho\rho_{1}\bar\eta, \ \ \ \ \mathcal{M}_{3}=\rho\rho_{2}\bar\eta,
\]
where $\rho=\max\{\bar c,1\}$ and $\rho_{1},\ \rho_{2}$ as in (\ref{e62}).
Letting $h\to 0$ and $\tau\to 0$, from (\ref{e63}) we obtain desired  estimate (\ref{et1}).  Therefore the theorem is proven.

Let us remark, the theorem \ref{t2} is also valid  for initial-boundary problems of the form (\ref{e1}) and (\ref{e2}) with nondegenerate in a domain of definition matrices  $A(x,t)$ and $B(x,t)$. To check this, it is enough to put $p=0$ and $l=0$.

\section{Numerical experiments}
\label{sect7}

For numerical solving of boundary problems  of the form (\ref{e1})-(\ref{ey1}) with mentioned above conditions, we have created a program in which the user can set input data, that is, matrix coefficients $A(x,t),\ B(x,t),\ C(x,t)$ and vector-functions  $f(x,t)$, $\phi(x)$ and $\psi(t)$; the domain of definition $U$; the values of the steps $h$ and $\tau$ and the orders of a spline  $m_{1}$ and $m_{2}$ with respect to each independent variable.

In this section we present the numerical results of solving of some problems of the form (\ref{e1})-(\ref{ey1}). Remark, that these examples were made only to demonstrate  the stability property of the difference scheme (\ref{e7}).

\begin{example}
Let us consider the system (\ref{e1}), in which
\[
A(x,t)=\left ( \begin{array}{cccccc}
1 & 0 & 0  & 1 & 0 & 0\\
0 & 0 & 0  & 0 & \exp(xt) & 0\\
0 & 0 & 0  & 0 & 0 & 0\\
1+xt & 0 & 0  & 0 & 0 & 0\\
0 & 0 & 1  & 0 & 0 & 0\\
0 & 0 & 0  & 0 & 0 & 0\\
\end{array}
\right ),\ \ \
B(x,t)=\left ( \begin{array}{cccccc}
0 & 0 & 0  & 0  & 0  & 0 \\
0 & 0 & 0  & 0 & \exp(\sin(\vartheta))  & 0 \\
0 & 0 & 0  & 0 & 0  & 1 \\
xt & 0 & 0  & 0 & 0  & 0 \\
0 & 0 & 0  & 0 & 0  & 0 \\
0 & 1 & 0  & 0 & 0  & 0 \\
\end{array}
\right ),
\]
\begin{equation}
C(x,t)=\left ( \begin{array}{cccccc}
0 & 0 & 0  & 0 & 0 & 0\\
0 & 0 & 0 & 0 & 2xt & 0\\
0 & 0 & 0  & 0 & 0 & \vartheta\\
1 & 0 & 0  & 0 & 0 & 0\\
0 & 0 & 1  & 0 & 0 & 0\\
0 & 0 & 0  & 0 & 0 & 0\\
\end{array}
\right ), \ \
f(x,t)=\left ( \begin{array}{l}
\exp(\vartheta) \\
\exp(xt)(x\exp(xt)+t\sin(\vartheta)+2xt)\\
\vartheta\\
x(1+t(\vartheta+1))\\
x\\
1\\
\end{array}
\right ).
\label{e64}
\end{equation}
with $\vartheta=x+t$. We know that the exact solution of our system with the data (\ref{e64}) is
\[
u(x,t)=(x,\ \exp(xt),\ \exp(\vartheta),\ xt,\ 1,\ \exp(xt)+\vartheta)^{\top}.
\]
With the help of nondegenerate in compact domain $\tilde U\subset\{(x,t)\in{\mathbb R}^{2},\ xt\neq -1\}$ matrices
\begin{equation}
P(x,t)=\left ( \begin{array}{cccccc}
0 & 1/\exp(xt) & 0  & 0 & 0 & 0\\
0 & 0 & 0  & 1/(1+xt) & 0 & 0\\
0 & 0 & 1  & 0 & 0 & 0\\
0 & 0 & 0  & 0 & 0 & 1\\
0 & 0 & 0  & 0 & 1 & 0\\
1 & 0 & 0  & 0 & 0 & 0\\
\end{array}
\right ),\ \
Q=\left ( \begin{array}{cccccc}
0 & 1 & 0  & 0 & 0 & 0\\
0 & 0 & 0  & 1 & 0 & 0\\
0 & 0 & 0  & 0 & 1 & 0\\
0 & -1 & 0  & 0 & 0 & 1\\
1 & 0 & 0  & 0 & 0 & 0\\
0 & 0 & 1  & 0 & 0 & 0\\
\end{array}
\right ),
\label{e65}
\end{equation}
the system (\ref{e1}) with the data (\ref{e64}) transforms to the canonical form  (\ref{ee1}), in which  $J_{1}(x,t)=\exp(\sin(\vartheta))/\exp(xt)$, $J_{2}(x,t)=xt/(1+xt)$ and  $M=N=\mathcal{O}_{2}$. To check this, it is enough to multiply the system (\ref{e1}) on the left by the matrix $P(x,t)$ and to make a change of variable: $v(x,t)=Qu(x,t)$. Remark, that in arbitrary compact domain $\tilde U$  the first condition of the theorem \ref{tp} is not fulfilled. To check this, it is enough to write down the characteristic equation of the system (\ref{e1}) with the data (\ref{e64}) and to see that its roots have common values in $\tilde U$.
Nevertheless, the system (\ref{e1}) with the data (\ref{e64}) can be cast to the form (\ref{ee1}) and therefore the conditions of the theorem \ref{t2} are fulfilled. Thus, in this case  the difference scheme (\ref{e7}) is stable in the domain $\tilde U$.
We show the results of numerical solving in the table 1. As stability estimate we take in this case the value of absolute error of the solution $\Delta u=\|u_{i,j}-v_{i,j}\|_{C(U_{\Delta})}$.

{\bf Table 1.}\\
\begin{center}
\begin{tabular}{|c|c|c|c|c|c|c|c|c|c|}
\hline
N & $h$ & $\tau$ & $t_{0}$ & $T$ & $x_{0}$ & $X$ & $m_{1}$ & $m_{2}$ & $\Delta u$\\
\hline \hline
1 & $10^{-1}$ & $10^{-1}$ & 0 & 1 & 0 & 1 & 2 & 2 & $2.07\times 10^{-2}$\\
2 & $10^{-1}$ & $10^{-1}$ & 0 & 1 & 0 & 1 & 3 & 2 & $2.07\times 10^{-2}$ \\
4 & $10^{-1}$ & $10^{-1}$ & 0 & 1 & 0 & 1 & 2 & 3 & $1.96\times 10^{-3}$ \\
5 & $10^{-1}$ & $10^{-1}$ & 0 & 1 & 0 & 1 & 3 & 3 & $1.96\times 10^{-3}$ \\
6 & $10^{-1}$ & $10^{-1}$ & 0 & 1 & 0 & 1 & 4 & 3 & $1.96\times 10^{-3}$ \\
7 & $10^{-1}$ & $10^{-1}$ & 0 & 1 & 0 & 1 & 3 & 4 & $1.91\times 10^{-4}$ \\
8 & $10^{-1}$ & $10^{-1}$ & 0 & 1 & 0 & 1 & 4 & 4 & $1.91\times 10^{-4}$ \\
9 & $10^{-1}$ & $10^{-1}$ & 0 & 1 & 0 & 1 & 5 & 5 & $1.91\times 10^{-5}$ \\
10 & $10^{-1}$ & $10^{-1}$ & 0 & 1 & 0 & 1 & 6 & 6 & $1.91\times 10^{-6}$ \\
11 & $10^{-1}$ & $10^{-1}$ & 0 & 1 & 0 & 1 & 7 & 7 & $1.93\times 10^{-7}$ \\
12 & $10^{-2}$ & $10^{-1}$ & 0 & 1 & 0 & 1 & 2 & 2 & $2.07\times 10^{-2}$ \\
13 & $10^{-1}$ & $10^{-2}$ & 0 & 1 & 0 & 1 & 2 & 2 & $1.28\times 10^{-3}$ \\
14 & $10^{-2}$ & $10^{-2}$ & 0 & 1 & 0 & 1 & 2 & 2 & $1.96\times 10^{-4}$ \\
15 & $5\times 10^{-3}$ & $5\times 10^{-3}$ & 0 & 1 & 0 & 1 & 2 & 2 & $4.95\times 10^{-5}$ \\
16 & $10^{-1}$ & $10^{-1}$ & 0 & 1 & 0 & 2 & 2 & 2 & $7.14\times 10^{-2}$ \\
17 & $10^{-1}$ & $10^{-1}$ & 0 & 1 & 0 & 2 & 3 & 3 & $1.31\times 10^{-2}$ \\
18 & $10^{-1}$ & $10^{-1}$ & 0 & 1 & 0 & 2 & 4 & 4 & $2.52\times 10^{-3}$ \\
19 & $10^{-1}$ & $10^{-1}$ & 0 & 1 & 0 & 2 & 5 & 5 & $4.99\times 10^{-4}$ \\

\hline
\end{tabular}
\end{center}
\vskip 0.3cm
\par\noindent

From the table 1 we see, that to achieve more accuracy of numerical solving, in some cases, it is enough to increase a degree of the spline with respect to one variable. This is important, because a degree of the spline impact on the dimension of the difference scheme (\ref{e7}) and therefore on a speed of calculations. What we learn  from the test 17 is that with using third-degree spline the accuracy of the numerical solutions in the domain $U=[0;1]\times [0;2]$  is of the second degree. This is explained by the impact of a large Lipchitz constant.
\end{example}

\begin{example}
Let us consider the system (\ref{e1}) with multiple characteristic curves with the following data
\[
A(x,t)={\rm diag} \{E_{5},\ 0,\ 1 \},\ \ \
B(x,t)={\rm diag} \{J_{1}(x,t),\ J_{2}(x,t),\ 1,\ 0 \},
\]
\[
J_{1}(x,t)=
\left ( \begin{array}{ccc}
\exp(\vartheta) & 1 & 0 \\
0 & \exp(\vartheta) & 1\\
0 & 0 & \exp(\vartheta)\\
\end{array}
\right ),\
J_{2}(x,t)=
\left ( \begin{array}{cc}
1+t\exp(x) & 1 \\
 0 & 1+t\exp(x) \\
\end{array}
\right ),
\]
\begin{equation}
C(x,t)=\left ( \begin{array}{ccccccc}
x^{2}+t & 0 & 1  & 1+xt & -\exp(\vartheta)  & 0 & 0\\
0 & x^{2} & xt & 0 & 0 & 1 & \vartheta\\
1 & 0 & 0  & 0 & 1 & xt  & 0\\
0 & 0 & 1  & 0 & 0 & 0  & 0\\
0 & 0 & 0  & 0 & 0 & 0  & 0\\
x\exp(\vartheta) & 1 & x^{2}t  & 0 & \vartheta & 0  & 0\\
0 & 0 & \exp(\vartheta)  & 0 & 0 & 1  & 0\\
\end{array}
\right ), \ \
\label{e66}
\end{equation}
\[
f(x,t)=\left ( \begin{array}{l}
\exp(2\vartheta)+(x^{2}+t)\exp(\vartheta)+2xt+(1+xt)(x-t)+1 \\
\exp(\vartheta)+2t+x^{3}+x^{2}t+3x^{2}t^{2}+x\exp(t)+x^{3}t+1\\
(2t+1)\exp(\vartheta)+2x+x^{2}t\exp(t)+1\\
2xt+t\exp(x)\\
0\\
x\exp(2\vartheta)+\exp(t)+2\vartheta+2x^{3}t^{2}\\
x^{2}+2xt\exp(\vartheta)+x\exp(t)\\
\end{array}
\right ).
\]
Here $\vartheta$ is defined in example  1. The exact solution of (\ref{e1})  with the data (\ref{e66}) is
\[
u(x,t)=(\exp(\vartheta), \ \vartheta,\ 2xt,\ x-t,\ 1,\ x\exp(t),\ x^{2}t)^{\top}.
\]
This system is given in canonical form. It has multiple nontrivial characteristic curves
\[
\lambda_{1,2,3}=-\exp(-\vartheta),\;\;\; \lambda_{4,5}=-1/(1+t\exp(x))
\]
and one zero characteristics  $\lambda_{6}=0$. In an arbitrary compact domain $\tilde U\subset U=\{(x,t),\ t\exp(x)>0\}$, the conditions of theorems \ref{tp} and \ref{t2} and remark \ref{z1} are valid. Therefore, the difference scheme (\ref{e7}) is stable in any domain $\tilde U$. The results of numerical solving are given in the table 2.

{\bf Table 2.}\\
\begin{center}
\begin{tabular}{|c|c|c|c|c|c|c|c|c|c|}
\hline
N & $h$ & $\tau$ & $t_{0}$ & $T$ & $x_{0}$ & $X$ & $m_{1}$ & $m_{2}$ & $\Delta u$\\
\hline \hline
1 & $10^{-1}$ & $10^{-1}$ & 0 & 1 & 0 & 1 & 2 & 2 & $3.54\times 10^{-2}$\\
2 & $10^{-1}$ & $10^{-1}$ & 0 & 1 & 0 & 1 & 3 & 2 & $4.96\times 10^{-3}$ \\
3 & $10^{-1}$ & $10^{-1}$ & 0 & 1 & 0 & 1 & 2 & 3 & $3.04\times 10^{-2}$ \\
4 & $10^{-1}$ & $10^{-1}$ & 0 & 1 & 0 & 1 & 3 & 3 & $3.34\times 10^{-3}$ \\
5 & $10^{-1}$ & $10^{-1}$ & 0 & 1 & 0 & 1 & 4 & 3 & $4.78\times 10^{-4}$ \\
6 & $10^{-1}$ & $10^{-1}$ & 0 & 1 & 0 & 1 & 3 & 4 & $2.91\times 10^{-3}$ \\
7 & $10^{-1}$ & $10^{-1}$ & 0 & 1 & 0 & 1 & 4 & 4 & $3.23\times 10^{-4}$ \\
8 & $10^{-1}$ & $10^{-1}$ & 0 & 1 & 0 & 1 & 5 & 3 & $6.19\times 10^{-4}$ \\
9 & $10^{-1}$ & $10^{-1}$ & 0 & 1 & 0 & 1 & 5 & 4 & $4.76\times 10^{-5}$ \\
10 & $10^{-1}$ & $10^{-1}$ & 0 & 1 & 0 & 1 & 5 & 5 & $3.21\times 10^{-5}$ \\
11 & $10^{-2}$ & $10^{-2}$ & 0 & 1 & 0 & 1 & 2 & 2 & $3.23\times 10^{-4}$ \\
12 & $10^{-1}$ & $10^{-1}$ & 0 & 2 & 0 & 2 & 2 & 2 & 7.97 \\
13 & $10^{-1}$ & $10^{-1}$ & 0 & 2 & 0 & 2 & 3 & 3 & $7.64\times 10^{-1}$ \\
14 & $10^{-1}$ & $10^{-1}$ & 0 & 2 & 0 & 2 & 4 & 4 & $7.46\times 10^{-2}$ \\
\hline
\end{tabular}
\end{center}
\vskip 0.3cm
\par\noindent
Remark, that in the tests 12, 13 and 14 on the error has influenced a large Lipchitz constant.
\end{example}

\section{Conclusion}

Numerous actual calculations show that spline collocation difference scheme given in this paper is a quite effective and gives sufficient accuracy. But, as was mentioned above, the difference scheme (\ref{e7}) in fact being a linear system of algebraic equations has the order $m_{1}m_{2}n$. Therefore increase of the degree of a used spline yields an increase of an order of this  system. Since our system by assumption has general form then  to solve it one needs to use, generally speaking, universal methods and therefore procedure of numerical solving might be quite laborious. In this case there is a need of a parallelization of numerical calculations. Following remark is in order.
The point is that many actual problems of mathematical physics and mechanics are described by  systems of  partial differential algebraic equations,  involving more than two independent variables. Examples are given by Sobolev's system, linearized Navier-Stokes system etc. \cite{dem}. Therefore the method given in this work  could be more marketable if one generalizes it to the case of many independent variables and this is the problem to be addressed.

\section*{Acknowledgements}
I would like to thank Svinin A.K. for valuable comments.



\begin{thebibliography}{99}


\bibitem{dem}
G. V. Demidenko, S. V. Uspenskii,
Equations and systems unsolved for the highest derivative, Nauchnaya kniga, Novosibirsk, 1998. (in Russian)


\bibitem{rush}
V. M. Rushchinskii,
Three-dimensional linear and nonlinear models of boiler-generators, Issues of identification and modeling,  (1968) 8-15. (in Russian)


\bibitem{cam}
S.L. Campbell, W. Marszalek,
The index of an infinite dimensional impliscit system,
Math. Comp. Model. Syst.,
5 (1999) 18-42.

\bibitem{luch}
W. Lucht, K. Strehmel, C. Eichler-Liebenow,
Indexes and special discretization methods for linear partial differential algebraic equations,
BIT, 39 (1999) 484-512.

\bibitem{deb}
K. Debrabant, K. Strehmel,
Convergence of Runge-Kutta methods applied to linear partial differential-algebraic equations,
Applied Numerical Mathematics 53 (2005) 213-229.

\bibitem{tisch}
C. Tischendorf,
Modeling circuit systems coupled with distributed semiconductor equations,
in:  Modeling, Simulation and Optimization of Integrated Circuits, International Series of Numerical Mathematics, Birkhauser,
Basel, 2003, Vol. 146, pp. 229-247.


\bibitem{tisch1}
C. Tischendorf,
Numerical analysis of DAEs from coupled circuit and semiconductor simulation,
Applied Numerical Mathematics, 53 (2005) 471-488.

\bibitem{han}
 M. Str\"omgren, M. Hanke,
On the numerical approximation of a degenerated hyperbolic system,
Mathematics and Computers in Simulation, 79 (2009) 1585-1602.
hyperbolic system

\bibitem{luch1}
W. Lucht,
Partial differential-algebraic systems of second order with symmetric convection,
Applied Numerical Mathematics, 53 (2005) 357-371.


\bibitem{gai1}
S. V. Gaidomak,  V. F. Chistyakov,
On Systems Other Than Cauchy-Kovalevskaya-Type Systems of Index $(1, k)$,
Vychisl. Tekhnol. 10 (2005) 45-59. (in Russian)



\bibitem{gai2}
S.V. Gaidomak,
Three-layer finite difference method for solving linear differential algebraic systems of partial differential equations,
Computational Mathematics and Mathematical Physics, 49 (2009) 1521-1534.


\bibitem{gai3}
S.V. Gaidomak,
Stability of an implicit difference scheme for a linear differential-algebraic system of partial differential equations,
~Computational Mathematics and Mathematical Physics, 50 (2010) 673-683.


\bibitem{gai4}
S.V. Gaidomak,
Spline collocation method for linear singular hyperbolic systems,
Computational Mathematics and Mathematical Physics, 48 (2008) 1161-1180.


\bibitem{chistbtv}
O.V. Bormotova, V.F. Chistyakov,
Numerical methods for solving and analyzing non Cauchy-Kovalevskaya-type systems,
Computational Mathematics and Mathematical Physics, 44 (2004) 1306-1313.

\bibitem{chist}
V. F. Chistyakov,
Differential-algebraic operators with a finite-dimensional kernel, Nauka, Novosibirsk, 1996. (in Russian)


\bibitem{gai}
S.V. Gaidomak,
The canonical structure of a pencil of degenerate matrix functions,
Russian Mathematics 56 (2012) 19-28.

\bibitem{zav}
Yu.S. Zavyalov, B.I. Kvasov, V.L. Miroshnichenko,
Spline function methods, Fizmatlit, Moscow, 1980, (in Russian)

\bibitem{sam}
A. A. Samarskii, A. V. Gulin,
Stability of difference schemes, Editorial URSS, Moscow, 2005. (in Russian)


\bibitem{ber}
M.V. Berezin,  N.P. Zhidkov,
Computation Methods, V.1, Nauka, Moscow , 1966. (in Russian)

\bibitem{gai5}
S.V. Gaidomak,
On some class of implisit spline collocation difference  schemes,
in: VI International conference ``Complex analysis and differential equations'', Ufa,
Institute of Mathematics with Computer Center of the Ufa Science Center of the Russian Academy of Sciences, 2011, pp. 50-51. (in Russian)


\bibitem{lank}
P. Lancaster,
Theory of Matrices, Academic Press, New York-London, 1969.

\bibitem{verb}
B.V. Verbitskii,
Some global property of matrix-functions that depend on several variables,
Russian Math. 22 (1978) 5-12.

\bibitem{gant}
F.R. Gantmacher,
The Theory of Matrices,
Chelsea, New York, 1959.

\bibitem{gai6}
S.V. Gaidomak,
On the numerical solution of a quasilinear algebraic-differential system,
Differential Equations. 45 (2009) 249-256.


\bibitem{vorob}
A. A. Vorobyev, M. Yu. Romanova,
On the estimates of norms powers of matrices,
Vestnic VGU, Seriya: Fizika. Matematika, 2 (2007) 83-85. (in Russian)


\bibitem{bulg}
A. Ya. Bulgakov,
An effectively calculable parameter for the stability property of a system of linear differential equations with constant coefficients,
Siberian Mathematical Journal, 21 (1980) 339-347.

\bibitem{nechep}
Yu.M. Nechepurenko,
Bounds for the matrix exponential based on the Lyapunov equation and limits of the Hausdorff set,
Computational Mathematics and Mathematical Physics, 42 (2002) 125-134.




\end{thebibliography}
\end{document}